\documentclass[12pt]{amsart}
\usepackage[dvips]{graphics}
\usepackage{amsfonts}
\usepackage{amssymb}
\usepackage{xypic}
\usepackage{a4}

\DeclareMathAlphabet{\eusm}{U}{}{}{}  
\SetMathAlphabet\eusm{normal}{U}{eus}{m}{n}
\SetMathAlphabet\eusm{bold}{U}{eus}{b}{n}

\DeclareMathAlphabet{\eufrak}{U}{}{}{}  
\SetMathAlphabet\eufrak{normal}{U}{euf}{m}{n}
\SetMathAlphabet\eufrak{bold}{U}{euf}{b}{n}

\newtheorem{theorem}{Theorem}[section]
\newtheorem{proposition}[theorem]{Proposition}
\newtheorem{lemma}[theorem]{Lemma}
\newtheorem{corollary}[theorem]{Corollary}

\theoremstyle{definition}
\newtheorem{definition}[theorem]{Definition}

\theoremstyle{remark}
\newtheorem{remark}[theorem]{Remark}

\numberwithin{equation}{section}

\newcommand{\alg}{\mathsf{A}}
\newcommand{\QG}{\mathbb{G}}

\newcommand{\Pol}{{\rm
  Pol}(\mathbb{G})}

\begin{document}

\title[Idempotent states on $U_q(2)$, $SU_q(2)$, and $SO_q(3)$]{Idempotent states on compact quantum groups and their classification on $U_q(2)$, $SU_q(2)$, and $SO_q(3)$}
\author{Uwe Franz}
\address{D\'epartement de math\'ematiques de Besan\c{c}on,
Universit\'e de Franche-Comt\'e 16, route de Gray, 25 030
Besan\c{c}on cedex, France}
\urladdr{http://www-math.univ-fcomte.fr/pp\underline{ }Annu/UFRANZ/}
\thanks{U.F.\ was supported by a Marie Curie Outgoing International
Fellowship of the EU (Contract Q-MALL MOIF-CT-2006-022137), an ANR Project
(Number ANR-06-BLAN-0015), and a Polonium cooperation}

\author{Adam Skalski}
\address{Mathematical Institute of the Polish Academy of Sciences,
ul.\'Sniadeckich 8, 00-956 Warszawa, Poland}
\email{a.skalski@impan.pl}

\author{Reiji Tomatsu}
\address{Department of Mathematics
Tokyo University of Science,
Yamazaki 2641, Noda, Chiba, 278-8510,
Japan}
\email{tomatsu\_reiji@mark\_ma.noda.tus.ac.jp}

\begin{abstract} Unlike for locally compact groups, idempotent states on locally compact quantum groups do not necessarily arise as Haar states of compact
quantum subgroups. We give a simple characterisation of those idempotent states on compact quantum groups which do arise as Haar
states on quantum subgroups. We also show that all idempotent states on the quantum groups  $U_q(2)$, $SU_q(2)$, and $SO_q(3)$
($q \in (-1,0)\cup (0,1]$) arise in this manner and list the idempotent states on the compact quantum semigroups $U_0(2)$,
$SU_0(2)$, and $SO_0(3)$. In the Appendix we provide a short new proof of coamenability of the deformations of classical compact
Lie groups based on their representation theory.
\end{abstract} 
 \maketitle

\section{Introduction}

It is well known that if $X$ is a locally compact topological semigroup, then the space of regular probability measures on $X$ possesses a natural convolution product. Analogously if $\alg$ is the algebra of continuous functions on a compact quantum semigroup, i.e.\ a unital $C^*$-algebra together with a coproduct (i.e., coassociative unital $*$-homomorphism) $ \Delta: \alg \to \alg \otimes \alg$ then one can consider a natural associative convolution product on the state space of $\alg$,
\[ \lambda \star \mu = (\lambda \otimes \mu)\circ \Delta, \;\;\; \lambda, \mu \in S(\alg).\]
It is natural to ask whether one can characterise the states which satisfy the idempotent property
\[ \mu \star \mu = \mu.\]
A particular and most important example of an idempotent state is the Haar state on a given compact quantum group 
in the sense of Woronowicz \cite{woronowicz98}. More general idempotent states arise naturally in considerations 
of C\'esaro limits of convolution operators on the algebras of functions on compact quantum groups, cf.\ 
\cite{franz+skalski08}. They are also an important ingredient in the construction of quantum hypergroups 
\cite{chapovsky+vainerman99} and occur as initial value $\varphi_0$ of convolution semigroups $(\varphi_t)_{t\ge 
0}$ of states on quantum groups, if one relaxes the initial condition $\varphi_0=\varepsilon$, cf.\ 
\cite{franz+schuermann00}. In the finite-dimensional case idempotent states are in 1-1 correspondence with 
quantum pre-subgroups studied in \cite{baaj+blanchard+skandalis99}. They may also be used in the study of Hopf 
images of representations of the algebras of functions on compact quantum groups introduced in 
\cite{Banica+Bichon10}.

For classical compact groups, Kawada and It\^o have proven that all idempotent measures are induced by Haar measures of compact
subgroups, see \cite[theorem 3]{kawada+ito40}. Later this result was extended to arbitrary locally compact topological groups,
see \cite{heyer77} and references therein. In \cite{pal96} Pal showed that this characterisation does not extend to quantum
groups by giving an example of an idempotent state on the Kac-Paljutkin quantum group which cannot arise as the Haar state on a
quantum subgroup. In \cite{franz+skalski08a} the first two authors began a systematic study of idempotent states on compact
quantum groups, later continued in \cite{franz+skalski09}. In particular we exhibited further examples of idempotent states on quantum groups that are not induced by Haar
states of quantum subgroups. Such examples are actually quite frequent and may be thought of as the dual manifestation of the
existence of non-normal subgroups. We also gave a characterisation of idempotent states on finite quantum groups in terms of sub
quantum hypergroups and established a 1-1 correspondence between idempotent states and quantum pre-subgroups of \cite{baaj+blanchard+skandalis99}.

In this work we continue the analysis began in \cite{franz+skalski08a} and show a simple and easy to check in concrete examples characterisation of those idempotent
states which arise as Haar states on quantum subgroups, so-called Haar idempotents. This, as in the finite-dimensional case leads to natural examples of Haar idempotents on group $C^*$-algebras of amenable discrete groups having non-normal subgroups. The second main result of the present paper is the
classification of all idempotent states on the compact quantum groups $U_q(2)$, $SU_q(2)$, and $SO_q(3)$ for $q \in
(-1,1]\setminus \{0\}$. Surprisingly, it turns out that they are all induced by quantum subgroups, cf.\ Theorems \ref{theo-uq2},
\ref{theo-suq2}, and \ref{theo-soq3}. As a byproduct we obtain the classification of quantum subgroups of the afore-mentioned
quantum groups. In the case of $SU_q(2)$ and $SO_q(3)$ these have been known thanks to the work of  of Podle\'s \cite{podles95}, but our approach provides a different proof.

For the value $q=0$, the quantum cancellation properties fail and $U_0(2)$, $SU_0(2)$, and $SO_0(3)$ are no 
longer compact quantum groups. But they can still be considered as compact quantum semigroups so that as 
explained above their state spaces have natural convolution products. \ Using this we determine all idempotent 
states on $U_0(2)$, $SU_0(2)$, and $SO_0(3)$, see Theorems \ref{theo-u02} and \ref{theo-su02} and 
\ref{theo-so03}. It turns out that in spite of the fact that the underlying $C^*$-algebras are isomorphic for all 
$q \in(-1,1)$, in the case $q=0$ there exist additional families of idempotent states, which do not appear when 
$q\not=0$. Interestingly, for $q=-1$ the Kac algebra $U_{-1}(2)$ admits non-Haar idempotent states. The full 
classification is however in this case more complicated and requires different techniques, so we refer it to the 
future work \cite{franz+skalski+tomatsu10}.

The detailed plan of the paper is as follows: in Section 2 we list the background results and definitions we need 
in the rest of the paper. In particular we recall the definitions of the quantum groups being the subject of the 
paper, discuss their representation theory and present explicit formulas for their Haar states. Although most of 
the facts presented in this section are well-known, the details of the respective representation theories are 
somewhat scattered through literature (and in the case of $U_q(2)$ even not fully recorded), so we decided to 
describe them here in detail to facilitate a coherent use of terminology and notation in the rest of the paper.  
Section 3 introduces idempotent states on compact quantum groups, provides the characterisation of those 
idempotent states which arise as  Haar states on quantum subgroups (extending the results for finite quantum 
groups given in \cite{franz+skalski08a}) and briefly discusses the cocommutative situation under the 
coamenability assumption. Section 4 contains main technical arguments of the paper and ends with the 
characterisation of all idempotent states on $U_q(2)$ for $q \in (-1,1]\setminus\{0\}$. In Section 5 we show how 
one can deduce the corresponding statements for $SU_q(2)$ and $SO_q(3)$. In view of Theorem 4.1 of 
\cite{franz+skalski09} this allows us immediately to conclude which of the coidalgebras of $C(SU_q(2))$ studied 
by the third named author in \cite{tomatsu08} are \emph{expected}, i.e.\ are images of a Haar state preserving 
conditional expectation  (we refer to both papers mentioned above for precise terminology).  Section 6 contains 
the classification of the idempotent states on compact quantum semigroups $U_0(2)$, $SU_0(2)$, and $SO_0(3)$. 
Finally in Section 7, we use the result of \cite{franz+skalski08a} showing that idempotent states on compact 
quantum groups are group-like projections in the dual quantum group, giving rise to algebraic quantum hypergroups 
by the construction given in \cite{landstad+vandaele07},  and discuss the quantum hypergroups associated to these 
group-like projections for the case of $SU_q(2)$. The appendix contains a new short proof of coamenability of 
deformations of classical compact Lie groups based on the representation theory developed in 
\cite{korogodski+soibelman98}.

\section{Notation and preliminaries}

The symbol $\otimes$ will denote the spatial tensor product of $C^*$-algebras,  $\odot$ will be reserved for the
purely algebraic tensor product. We write $\mathbb{N}_0$ or $\mathbb{Z}_+$ for $\mathbb{N}\cup \{0\}$.

\subsection{Compact quantum groups}

The notion of compact quantum groups has been introduced in \cite{woronowicz87}. Here we adopt a rewording of the 
definition from \cite{woronowicz98} (Definition 2.1 of that paper).

\begin{definition}
A \emph{$C^*$-bialgebra} (the algebra of continuous functions on a compact quantum semigroup) is a pair $(\mathsf{A}, \Delta)$, where
$\mathsf{A}$ is a unital $C^*$-algebra,
 $\Delta:\mathsf{A} \to \mathsf{A} \otimes \mathsf{A}$ is a unital,
 $*$-homomorphic map which is coassociative
\[ (\Delta \otimes \textup{id}_{\mathsf{A}})\circ \Delta = (\textup{id}_{\mathsf{A}} \otimes \Delta)
\circ\Delta.
\]
If the quantum cancellation properties
\[
\overline{{\rm Lin}}((1\otimes \mathsf{A})\Delta(\mathsf{A}) ) = \overline{{\rm Lin}}((\mathsf{A} \otimes 1)\Delta(\mathsf{A}) )
= \mathsf{A} \otimes \mathsf{A},
\]
are satisfied, then the pair $(\mathsf{A}, \Delta)$ is called the algebra of continuous functions on a \emph{compact
  quantum group}.
\end{definition}

In quantum group theory it is quite common to write
$\mathsf{A}=C(\mathbb{G})$, where $\mathbb{G}$ denotes the underlying compact quantum group --- note
however that the symbol $\mathbb{G}$ is only defined indirectly.

The map $\Delta$ is called the coproduct of $\mathsf{A}$, it induces the convolution product
\[\lambda\star \mu:=(\lambda\otimes \mu)\circ\Delta, \;\;\; \lambda,\mu\in\mathsf{A}^*.\]

If $\QG$ is a compact quantum group then a unitary $U \in M_n(C(\QG))$ is called a {\emph{finite-dimensional 
unitary representation of}} $\mathbb{G}$ if for all $i,j=1,\ldots,n$ we have $\Delta(U_{ij}) = \sum_{k=1}^n 
U_{ik} \otimes U_{kj}$.  It is said to be \emph{irreducible}, if the only matrices $T\in M_n(\mathbb{C})$ with 
$TU=UT$ are multiples of the identity matrix.

Possibly the most important feature of compact quantum groups is the existence of the dense $*$-subalgebra $\mathcal{A}\subset C(\QG)$
(the algebra of matrix coefficients of irreducible unitary representations of $\QG$), which is in fact a Hopf
$*$-algebra - so for example $\Delta: \mathcal{A} \to \mathcal{A} \odot
\mathcal{A}$. This $*$-Hopf algebra is also denoted by $\mathcal{A}=\Pol$, and treated as the analog of polynomial functions of
$\mathbb{G}$.

Another fact of the crucial importance is given in the following result, Theorem 2.3 of
\cite{woronowicz98}.

\begin{theorem} \label{prop-haar}
Let $\QG$ be a compact quantum group and let $\alg=C(\QG)$. There exists a unique state $h \in \alg^*$ (called the \emph{Haar state} of
$\QG$) such that for all $a \in \mathsf{A}$
\[
(h \otimes \textup{id}_{\mathsf{A}})\circ  \Delta (a) = ( \textup{id}_{\mathsf{A}} \otimes h)\circ  \Delta (a) = h(a) 1.
\]
\end{theorem}

 The algebra $C(\QG)$ is said to be \emph{in reduced form} if the Haar state $h$ is faithful. If it is not the case we can
always quotient out the kernel of $h$. This procedure in particular does not influence the underlying Hopf 
$*$-algebra $\mathcal{A}$; in fact the reduced object may be viewed as the natural completion of $\mathcal{A}$ in 
the GNS representation with respect to $h$ (as opposed for example to the universal completion of $\mathcal{A}$, 
for details see \cite{bedos+murphy+tuset01}). In general the reduced and universal object need not coincide. This 
leads to certain technical complications which are not of essential importance in our context (for example if a 
discrete group $\Gamma$ is not amenable, the reduced $C^*$-algebra of $\Gamma$ is a proper quantum subgroup of 
the universal $C^*$-algebra of $\Gamma$, even though they have `identical' Haar states). To avoid such 
difficulties we focus on the class of coamenable compact quantum groups (\cite{bedos+murphy+tuset01}, see the 
Appendix to this paper for more information), for which the reduced and universal $C^*$-algebraic completions of 
$\mathcal{A}$ are naturally isomorphic. All deformations of classical compact Lie groups, so in particular 
quantum groups $U_q(2)$, $SU_q(2)$ and $SO_q(3)$ we consider in Sections 4 and 5 are known to be coamenable 
(\cite{Banica99}). Banica's proof is based on the fusion rules appearing in the quantum group theoretic 
representation theory; in the Appendix we give an alternative proof based on the $C^*$-algebraic representation 
theory developed in \cite{korogodski+soibelman98}.

The following definition was introduced by Podle\'s in the context of compact matrix pseudogroups (Definition 1.3 of
\cite{podles95}).

\begin{definition} A compact quantum group $\QG'$ is said to be a quantum subgroup of
a compact quantum group $\QG$ if there exists a surjective compact quantum group morphism  $j:C(\QG) \to C(\QG')$, i.e.\ a surjective
unital
$*$-homomorphism $j:C(\QG) \to C(\QG') $ such that
\begin{equation} \Delta_{C(\QG')}\circ j = (j \otimes j) \circ\Delta_{C(\QG)}.\label{comor}\end{equation}
\end{definition}

Strictly speaking, one should consider the pairs $(\QG',j)$, since $C(\QG)$ can contain several copies of 
$C(\QG')$ with different morphisms (in the same way as a classical group can have different, but isomorphic 
subgroups). We will not distinguish between $(\QG_1,j_1)$ and $(\QG_2,j_2)$ if there exists an isomorphism of 
quantum groups $\Theta:C(\QG_1)\to C(\QG_2)$ such that $\Theta\circ j_1=j_2$. Note that such isomorphic pairs 
induce the same idempotent state $\phi=h_{\QG_1}\circ j_1=h_{\QG_2}\circ j_2$, since uniqueness of the Haar 
states implies $h_{\QG_1}=h_{\QG_2}\circ \Theta$.

Any coamenable compact quantum group contains itself and the trivial compact quantum group $\{e\}$ as quantum subgroups;
further these two quantum subgroups will be called trivial. If $G$ is a compact group and a compact quantum group $\QG$
contains $G$ as a quantum subgroup, via a morphism $j:C(\QG)\to C(G)$, we will  simply say that $G$ is a subgroup of
$\QG$.

\begin{definition}
A quantum subgroup $(\QG',j)$ of a compact quantum group $\QG$ is called normal if the images of the conditional 
expectations
\begin{eqnarray*}
E_{\QG/\QG'} &=& \big({\rm id}\otimes (h_{\QG'}\circ j)\big)\circ \Delta, \\
E_{\QG'\backslash \QG} &=& \big((h_{\QG'}\circ j)\otimes {\rm id}\big)\circ \Delta,
\end{eqnarray*}
coincide, cf.\ \cite[Proposition 2.1 and Definition 2.2]{wang08}. In this case
the quotient $C(\QG/\QG'):=E_{\QG/\QG'}(C(\QG))$ has a natural structure of the algebra of continuous functions on a compact quantum group. \label{defnormal}
\end{definition}

\subsection{$q$-Numbers}
Let $q\not=1$. We will use the following notation for $q$-numbers,
\begin{eqnarray*}
(x;q)_n &=& (1-x)(1-qx)\cdots(1-q^{n-1}x), \\
\left[\begin{array}{cc} n \\ k \end{array}\right]_q &=& \frac{(q;q)_n}{(q;q)_k(q;q)_{n-k}},
\end{eqnarray*}
for $n\in\mathbb{Z}_+$, $0\le k \le n$, $x \in \mathbb{R}$.

\subsection{The Woronowicz quantum group $SU_q(2)$}

\cite{woronowicz87,woronowicz87b} For $q\in\mathbb{R}$, we denote by ${\rm
  Pol} (SU_q(2))$ the $*$-bialgebra generated by $\alpha$ and $\gamma$, with the relations
\begin{gather*}
\alpha \gamma = q \gamma \alpha, \qquad \alpha \gamma^* = q \gamma^*
\alpha, \qquad\gamma^* \gamma = \gamma \gamma^*, \\
\gamma^* \gamma + \alpha^*
\alpha = 1, \qquad \alpha \alpha^* - \alpha^*\alpha = (1-q^2) \gamma^* \gamma,
\end{gather*}
and comultiplication and counit defined by setting
\[
\Delta (\alpha) = \alpha\otimes\alpha - q \gamma^*\otimes \gamma, \quad
\Delta(\gamma) = \gamma\otimes \alpha + \alpha^* \otimes \gamma,
\]
and $\varepsilon(\alpha) = 1$, $\varepsilon(\gamma) =0$. For $q\not=0$, ${\rm
  Pol} (SU_{q}(2))$ admits an antipode. On the generators, it acts as
\[
S(\alpha)=\alpha^* \qquad\mbox{ and }\qquad S(\gamma)=-q\gamma.
\]

Denote the universal enveloping $C^*$-algebra of ${\rm Pol}(SU_q(2))$ by
$C(SU_q(2))$, then $\Delta$ extends uniquely to a non-degenerate coassociative
homomorphism $\Delta:C(SU_q(2))\to C(SU_q(2))\otimes C(SU_q(2))$, and the pair
$(C(SU_q(2)),\Delta)$ is a $C^*$-bialgebra. For $q\not=0$, $SU_q(2)$ is even a compact quantum group.

Note that the mapping $\alpha\mapsto \alpha^*$ and $\gamma\mapsto q\gamma^*$
induces isomorphisms ${\rm Pol}(SU_{1/q}(2))\to{\rm Pol} (SU_q(2))$,
$C(SU_{1/q}(2))\to C(SU_q(2))$, therefore it is sufficient to consider $q\in[-1,1]$.

\subsubsection{Representation theory of $C(SU_q(2))$} \cite{woronowicz87b,
vaksman+soibelman88} The $C^*$-algebra $C(SU_q(2))$ has two families of irreducible
representations. The first family consists of the one-dimensional representations $\rho_\theta$,
$0\le \theta < 2\pi$, given by
\begin{eqnarray*}
\rho_\theta(\alpha) &=& e^{i\theta}, \\
\rho_\theta(\gamma) &=& 0.
\end{eqnarray*}

The other family consist of infinite-dimensional representations $\pi_\theta$, $0\le \theta< 2\pi$, acting on a separable Hilbert space $\eufrak{h}$ by
\begin{eqnarray*}
\pi_\theta(\alpha)e_n &=& \left\{
\begin{array}{ccl}
\sqrt{1-q^{2n}}\,e_{n-1} & \mbox{ if } & n>0, \\
0 & \mbox{ if } & n=0,
\end{array}\right. \\
\pi_\theta(\gamma)e_n &=& e^{i\theta}q^ne_n,
\end{eqnarray*}
where $\{e_n;n\in\mathbb{N}_0\}$ is an orthonormal basis for $\eufrak{h}$.

This list is complete, i.e.\ any irreducible representation of $C(SU_q(2))$ is unitarily equivalent
to a representation in one of the two families above (Theorem 3.2 of \cite{vaksman+soibelman88}).
It is known that the $C^*$-algebra $C(SU_q(2))$ is of type I, therefore any representation can be
written as a direct integral over the irreducible representations given above.

\subsubsection{Quantum subgroups of $SU_q(2)$}

Let ${\rm Pol}(U(1))$ denote the $*$-algebra generated by one unitary $w$, $ww^*=w^*w=1$. With
$\Delta(w)=w\otimes w$, $\varepsilon(w)=1$, $S(w)=w^*$, this becomes a $*$-Hopf algebra. Its
enveloping $C^*$-algebra $C(U(1))$ is the algebra of continuous functions on a compact group $U(1)$. Note that the $*$-algebra
homomorphism ${\rm Pol}(SU_q(2)) \to {\rm Pol}(U(1))$ defined by $\alpha\mapsto w$, $\gamma\mapsto
0$ extends to a surjective compact quantum group morphism $j:C(SU_q(2))\to C(U(1))$, i.e.\ $U(1)$
is a quantum subgroup of $SU_q(2)$. Furthermore, Podle\'s in Theorem 2.1 of \cite{podles95} has
shown that, for $q\in(-1,0)\cup(0,1)$, $U(1)$ and its closed subgroups are the only non-trivial
quantum subgroups of $SU_q(2)$. This will also follow from the results in Section 4.

There exists a second morphism $j':C(SU_q(2))\to C(U(1))$, determined by $j':\alpha\mapsto w^*,
\gamma\mapsto 0$. But we do not need to distinguish the pairs $(C(U(1)),j)$ and $(C(U(1)),j')$,
since they are related by the automorphism $\Theta$ of $C(U(1))$ with $\Theta(w^k) = (w^*)^k$,
$\Theta\big((w^*)^k\big) = w^k$, $k \in \mathbb{N}_0$.

\subsubsection{Representations of $SU_q(2)$}

Let $q\in(-1,0)\cup(0,1)$. We recall a few basic facts about the representations of $SU_q(2)$, for more details 
see \cite{woronowicz87,woronowicz87b,vaksman+soibelman88,masuda+al88,koornwinder89}. For each non-negative 
half-integer $s\in\frac{1}{2}\mathbb{Z}_+$ there exists a $2s+1$-dimensional irreducible unitary representation 
$u^{(s)}=(u^{(s)}_{k\ell})_{-s\le k,\ell\le s}$ of $SU_q(2)$, which is unique up to unitary equivalence. Note 
that indices $k,\ell$ run over the set $\{-s,-s+1,\ldots,s-1,s\}$, they are integers if $s\in\mathbb{Z}_+$ is 
integer, and half-integer if $s\in(\frac{1}{2}\mathbb{Z}_+)\backslash\mathbb{Z}_+$ is half-integer. This 
convention is also used further in the paper.

The matrix coefficients $u^{(s)}_{k\ell}$, $s\in\frac{1}{2}\mathbb{Z}_+$, $-s\le k,\ell\le s$, span
${\rm Pol} (SU_q(2))$ and are linearly dense in $C(SU_q(2))$, therefore will be sufficient for
our calculations. We have
\begin{gather*}
u^{(0)}= (1), \qquad u^{(1/2)} = \left(\begin{array}{cc} \alpha & -q\gamma^*
    \\ \gamma & \alpha^* \end{array}\right), \\
\qquad u^{(1)} = \left(\begin{array}{ccc} \alpha^2 &
-q\sqrt{1+q^2}\gamma^*\alpha & q^2 (\gamma^*)^2 \\
\sqrt{1+q^2}\gamma\alpha & 1-(1+q^2) \gamma^*\gamma &
-q\sqrt{1+q^2}\alpha^*\gamma^* \\
\gamma^2 & \sqrt{1+q^2} \alpha^*\gamma & (\alpha^*)^2\end{array}\right),
\end{gather*}
and the matrix coefficients of the higher-dimensional representations are of the
form
\[
u^{(s)}_{k\ell}=\left\{\begin{array}{ccl}
\alpha^{-k-\ell}p^{(s)}_{k\ell} \gamma^{k-\ell} & \mbox{ for } & k+\ell\le
0,\, k\ge \ell, \\
\alpha^{-k-\ell}p^{(s)}_{k\ell} (\gamma^*)^{\ell-k} & \mbox{ for } & k+\ell\le
0,\, k\le \ell, \\
(\alpha^*)^{k+\ell}p^{(s)}_{k\ell} \gamma^{k-\ell} & \mbox{ for } & k+\ell\ge
0,\,k\ge \ell, \\
(\alpha^*)^{k+\ell}p^{(s)}_{k\ell} (\gamma^*)^{\ell-k} & \mbox{ for } & k+\ell\ge 0,\, k\le \ell,
\end{array}\right.
\]
where $p^{(s)}_{k\ell}$ is a polynomial in $\gamma^*\gamma$.

In particular, for $s$ an integer,
$u^{(s)}_{00}=p_{s}(\gamma^*\gamma;1,1;q^2)$ is the little $q$-Legendre polynomial
and $u^{(s)}_{s0}=\sqrt{\left[\begin{array}{c} 2s \\ s \end{array}\right]_{q^2}}(\alpha^*)^s\gamma^s$.

If we define a $\mathbb{Z}$-grading on ${\rm Pol} (SU_q(2))$ by
\[
\deg \alpha = \deg\alpha^*=0,\quad \deg \gamma=1,\quad
\deg\gamma^*=-1,
\]
then we have
\[
\deg u^{(s)}_{k\ell} = k-\ell.
\]
With this grading, it is straight-forward to verify the following formula for the square of the antipode on 
homogeneous elements,
\begin{equation}\label{S square}
S^2(a) = q^{2\,\deg a} a.
\end{equation}

\subsubsection{The Haar state of $SU_q(2)$}

As stated in Theorem \ref{prop-haar}, there exists a unique invariant state $h$ on the compact
quantum group $SU_q(2)$, called the Haar state. It was first computed by Woronowicz in \cite{woronowicz87}; it is the identity on the
one-dimensional representation, and vanishes on the matrix coefficients of all other irreducible
representations, i.e.\
\[
h\left(u^{(s)}_{k\ell}\right) = \delta_{0s}
\]
for $s\in\frac{1}{2}\mathbb{Z}_+$, $-s\le k,\ell\le s$. On polynomials $p(\gamma^*\gamma)\in\mathbb{C}[\gamma^*\gamma]$, it is
equal to Jackson's $q$-integral (\cite{koornwinder89}),
\[
h\big(p(\gamma^*\gamma)\big) = (1-q^2) \sum_{k=0}^\infty q^{2k} p(q^{2k}) =:
\int_0^1 p(x) {\rm d}_{q^2}x.
\]

\subsection{The compact quantum group $SO_q(3)$}
A compact quantum group $\QG'$ is called a quotient group of $\QG$, if there exists an injective morphism of 
quantum groups $j:C(\QG')\to C(\QG)$. The compact quantum group $SO_q(3)$ can be defined as the quotient of 
$SU_q(2)$ by the quantum subgroup $\mathbb{Z}_2$, cf.\ \cite{podles95}. ${\rm Pol} (SO_q(3))$ is the subalgebra 
of ${\rm Pol} (SU_q(2))$ spanned by the matrix coefficients of the unitary irreducible representations of 
$SU_q(2)$ with integer label, and $C(SO_q(3))$ its norm closure. The Haar state on $SO_q(3)$ is simply the 
restriction of the Haar state on $SU_q(2)$.

Podle\'s  \cite{podles95} has shown that $SO_q(3)$ and $SO_{-q}(3)$ are isomorphic.

\subsubsection{The semigroup case $q=0$} We define ${\rm Pol} (SO_0(3))$ as the
unital $*$-subalgebra of ${\rm Pol} (SU_0(2))$ generated by $\alpha^2$,
$\gamma^2$, $\gamma\alpha$, $\gamma^*\alpha$, and $\gamma^*\gamma$, i.e.\
\[
{\rm Pol} (SO_0(3))={\rm span}\,\big\{(\alpha^*)^r \gamma^k
\alpha^s,(\alpha^*)^r (\gamma^*)^k \alpha^s:
r,k,s\in\mathbb{Z}_+~\mbox{s.t.}~r+k+s\mbox{ even}\}.
\]
Since
\begin{eqnarray*}
\Delta\big((\alpha^*)^r \gamma^k \alpha^s \big) &=& \sum_{\kappa=0}^k
(\alpha^*)^{r+k-\kappa}\gamma^\kappa \alpha^s \otimes (\alpha^*)^r
\gamma^{k-\kappa}\alpha^{s+\kappa}, \\
\Delta\big((\alpha^*)^r (\gamma^*)^k \alpha^s \big) &=& \sum_{\kappa=0}^k
(\alpha^*)^r (\gamma^*)^\kappa \alpha^{s+k-\kappa} \otimes (\alpha^*)^{r+\kappa}
(\gamma^*)^{k-\kappa}\alpha^s,
\end{eqnarray*}
this is a sub $*$-Hopf algebra in ${\rm Pol} (SU_0(2))$. The $C^*$-bialgebra $C(SO_0(3))$ is then defined as the 
norm closure of ${\rm Pol}(SO_0(3))$ in $C(SU_0(2))$. Note that $C^*$-algebras $C(SU_q(2))$ are isomorphic for 
all $q\in (-1,1)$ \cite[Theorem A2.2]{woronowicz87b}.

\subsubsection{The conditional expectation $E:C(SU_q(2))\to C(SO_q(3))$}\label{subsub-cond-exp}
Looking at the defining relations of $SU_q(2)$, it is clear that
$\vartheta:\alpha\to-\alpha,\gamma\to-\gamma$ extends to a unique $*$-algebra
automorphism of $C(SU_q(2))$. Therefore $E=\frac{1}{2}({\rm id}+\vartheta)$
defines a completely positive unital map from $C(SU_q(2))$ to itself. If
$\phi_2=h_{\mathbb{Z}_2}\circ j$ denotes the idempotent state on $SU_q(2)$
induced by the Haar measure of $\mathbb{Z}_2$, with $j:SU_q(2) \to C(\mathbb{Z}_2)$ the corresponding surjective morphism (see Definitions \ref{defidemp} and \ref{defHaaridemp}), then we can write $E$ also as
\[
E=({\rm id}\otimes \phi_2)\circ \Delta.
\]
Checking
\[
E(u^{(s)}_{k\ell}) = \left\{
\begin{array}{cl}
u^{(s)}_{k\ell} & \mbox{ if } s\in\mathbb{Z}_+, \\
0 & \mbox{ else}.
\end{array}\right.
\]
we can show that the range of $E$ is equal to $C(SO_q(3))$. Furthermore, $E$
satisfies
\[
\Delta\circ E = ({\rm id}\otimes E)\circ\Delta = (E\otimes {\rm id})\circ
\Delta = (E\otimes E)\circ \Delta.
\]

\subsubsection{Quantum subgroups of $SO_q(3)$}
The restriction of the morphism $j:C(SU_q(2))\to C(U(1))$ to $C(SO_q(3))$ is no longer surjective, its range is 
equal to the subalgebra $\{f\in C(U(1)): \forall_{z\in U(1)}\, f(z)=f(-z)\}=C(U(1)/\mathbb{Z}_2)$. Since 
$U(1)/\mathbb{Z}_2\cong U(1)$, we see that $SO_q(3)$ contains $U(1)\cong SO(2)$ and its closed subgroups as 
quantum subgroups. Podle\'s  \cite{podles95} has shown that these are the only non-trivial quantum subgroups of 
$SO_q(3)$. Again this can be deduced from the results of Section 5.

\subsection{The compact quantum group $U_q(2)$}

\cite{koelink91,wysoczanski04,zhang+zhao05} Let $q\in \mathbb{R}$. Then ${\rm
Pol} (U_q(2))$ is defined as the $*$-bialgebra generated by $a$, $c$, and
$v$, with the relations
\begin{gather*}
av=va, \qquad cv = vc, \qquad cc^* = c^*c , \\
ac=qca, \qquad a c^* = qc^* a, \qquad vv^* = v^*v= 1, \\
aa^* + q^2cc^* =  1 = a^*a + c^*c, \\
\Delta(a) = a\otimes a-qc^*v^* \otimes c, \quad \Delta (c) = c\otimes a + a
\otimes c, \quad \Delta(v)=v\otimes v, \\
\varepsilon(a)=\varepsilon(v)=1, \qquad \varepsilon(c)=0.
\end{gather*}
For $q\not=0$, ${\rm Pol} (U_q(2))$ admits an antipode, given by
\[
S(a) = a^*, \qquad S(v)=v^*,\qquad S(c)= -qcv,
\]
on the generators.

Denote the universal enveloping $C^*$-algebra of ${\rm Pol}(U_q(2))$ by
$C(U_q(2))$, then $\Delta:{\rm Pol} (U_q(2))\to {\rm Pol}(U_q(2))\odot{\rm
  Pol}(U_q(2))$ extends uniquely to a non-degenerate coassociative
homomorphism $\Delta:C(U_q(2))\to C(U_q(2))\otimes C(U_q(2))$, and the pair
$(C(U_q(2)),\Delta)$ is a $C^*$-bialgebra. For $q\not=0$, $U_q(2)$ is even
a compact quantum group. It is again sufficient to consider $q\in [-1,1]$,
since $U_q(2)$ and $U_{1/q}(2)$ are isomorphic.

\subsubsection{Quantum subgroups of $U_q(2)$}
The mapping $a\mapsto\alpha$, $c\mapsto \gamma$, $v\mapsto 1$ extends to a surjective compact quantum group 
morphism $C(U_q(2))\to C(SU_q(2))$ and shows that $SU_q(2)$ is a quantum subgroup of $U_q(2)$. The $C^*$-algebra 
$C(U_q(2))$ is isomorphic to the tensor product of $C(SU_q(2))$ and $C(U(1))$. Moreover the compact quantum 
group, $U_q(2)$  is equal to a twisted product of $SU_q(2)$ and $U(1)$, cf.\ \cite{wysoczanski04}, written as 
$U_q(2)=SU_q(2)\ltimes_\sigma U(1)$.

Another quantum subgroup of
$U_q(2)$ is the two-dimensional torus. Denote by ${\rm Pol}(\mathbb{T}^2)$ the
$*$-Hopf algebra generated by two commuting unitaries, i.e.\ by $w_1,w_2$ with
the relations
\begin{gather*}
w_1w_1^*=1=w_1^*w_1, \quad w_2w_2^*=1=w_2^*w_2,\quad w_1w_2=w_2w_1,\quad
w_1w_2^*=w_2^*w_1, \\
\Delta (w_1)=w_1\otimes w_1,\quad \Delta(w_2)=w_2\otimes w_2, \quad \varepsilon(w_1)=\varepsilon(w_2)=.
\end{gather*}
Then $C(\mathbb{T}^2)$ is the
$C^*$-enveloping algebra of ${\rm Pol}(\mathbb{T}^2)$. Then the mapping $a\mapsto w_1$, $c\mapsto 0$,
$v\mapsto w_2$ extends to a unique surjective compact quantum group morphism
$C(U_q(2))\to C(\mathbb{T}^2)$.

We will see that the twisted products $SU_q(2)\ltimes_\sigma \mathbb{Z}_n$,
$n\in\mathbb{N}$, the
torus $\mathbb{T}^2$, and its closed subgroups are the only non-trivial quantum subgroups of $U_q(2)$, cf.\ Corollary \ref{uq2-subgroups}.

\subsubsection{Representations of $U_q(2)$}

Unitary irreducible representations of $U_q(2)$ can be obtained as tensor
product of unitary irreducible representations of $SU_q(2)$ with
representations of $U(1)$, cf.\ \cite{wysoczanski04}. In this way one obtains
the following family of unitary irreducible representations of $U_q(2)$,
\[
v^{(s,p)} = \left(u^{(s)}_{k\ell}v^{p+s+\ell}\right)_{-s\le k,\ell\le s}
\]
for $p\in\mathbb{Z}$, $s\in\frac{1}{2}\mathbb{Z}_+$. The matrix coefficients of
these representations clearly span ${\rm Pol} (U_q(2))$. Therefore they are
dense in $C(U_q(2))$ and will be sufficient for the calculations in this
paper.

Assume $q\not=0$. We want to compute the action of the square of the antipode on the matrix coefficients of the 
unitary irreducible representations defined above. Since we have $S^2(a)=a$, $S^2(c)=q^2c$, and $S^2(v)=v$, we get
\begin{equation}\label{u2-s-squared}
S^2(u^{(s)}_{k\ell}v^{p+s+\ell}) =q^{2(k-\ell)}u^{(s)}_{k\ell}v^{p+s+\ell}
\end{equation}
$p\in\mathbb{Z}$, $s\in\frac{1}{2}\mathbb{Z}_+$, and $-s\le k,\ell\le s$.

\subsubsection{The Haar state of $U_q(2)$}

The Haar state $h$ on $U_q(2)$ can be written as a tensor product of the Haar state on $SU_q(2)$ and the Haar 
state on $U(1)$, it acts on the matrix coefficients of the unitary irreducible representations given above as
\[
h\left(u^{(s)}_{k\ell}v^{p+s+\ell}\right) = \delta_{0s}\delta_{0p}
\]
for $s\in\frac{1}{2}\mathbb{Z}_+$, $-s\le k,\ell\le s$, $p\in\mathbb{Z}$.

\subsection{Multiplicative domain of a completely positive unital map}
The following result by Choi on multiplicative domains will be useful for us.
\begin{theorem}\label{mult-dom}\cite[Theorem 3.1]{choi74}, \cite[Theorem 3.18]{paulsen02}
Let $T:A\to B$ be a completely positive unital linear map between two $C^*$-algebras $A$ and $B$.
Set
\[
D_T=\{a\in A: T(aa^*)=T(a)T(a^*), T(a^*a)=T(a^*)T(a)\}.
\]
Then we have
\[
T(ab)=T(a)T(b) \qquad\mbox{ and }\qquad T(ba)=T(b)T(a)
\]
for all $a\in D_T$ and $b\in A$.
\end{theorem}

\section{Idempotent states on compact quantum groups}

In this section we formally introduce the notion of idempotent states on a $C^*$-bialgebra, provide a characterisation of those
idempotent states on compact quantum groups which arise as Haar states on quantum subgroups and discuss commutative and
cocommutative cases.

\begin{definition} Let $(\mathsf{A},\Delta)$ be a $C^*$-bialgebra. A state
  $\phi\in \mathsf{A}^*$ is called an {\em idempotent state} if
\[
(\phi \otimes \phi) \circ\Delta = \phi,
\]
i.e.\ if it is idempotent for the convolution product. \label{defidemp}
\end{definition}

If $\mathsf{A}=C(\QG)$ for a compact quantum group $\QG$, we will also simply say that $\phi$ as above is an idempotent state on $\QG$.

\begin{definition} \label{defHaaridemp}
A state $\phi\in C(\QG)^*$ is called a {\em Haar state on a quantum subgroup of
$\QG$} (or a {\em Haar idempotent}) if there exists a quantum subgroup
($\QG',j$) of $\QG$ and $\phi = h_{\QG'} \circ j$, where  $h_{\QG'}$ denotes the Haar state on $\QG'$.
\end{definition}

It is easy to check that each Haar state on a quantum subgroup of $\QG$ is idempotent. It
follows from the example of Pal in \cite{pal96} and our work in \cite{franz+skalski08a} that not
every idempotent state is a Haar idempotent. We have the following simple characterisation,
extending Theorem 4.5 of \cite{franz+skalski08a}.

\begin{theorem}
Let $\QG$ be a compact quantum group, let $\phi\in C(\QG)^*$ be an idempotent state and let $N_{\phi} = \{a \in C(\QG) : \phi (a^*a)
=0\}$ denote the null space of $\phi$. Then $\phi$ is a Haar idempotent if and only if $N_{\phi}$ is a two-sided (equivalently,
selfadjoint) ideal.
\end{theorem}

\begin{proof}
It is an easy consequence of the Cauchy-Schwarz inequality that $N_{\phi}$ is a left ideal; thus it is a two-sided ideal if and
only if it is selfadjoint.

Write $\mathsf{A}=C(\QG)$.
Suppose first that $\phi$ is a Haar idempotent, i.e.\ there exists a compact quantum group $\QG'$ and a surjective compact
quantum group morphism $j:\mathsf{A} \to C(\QG')$ such that $\phi = h_{\QG'} \circ j$. Recall that we assumed $h_{\QG'}$ to be faithful,
so that $N_{\phi} = \{a \in \mathsf{A}: j(a^*a)=0\}=\{a\in \mathsf{A}: j(a) =0\}$, which is obviously self-adjoint.

Suppose then that $N_{\phi}$ is a two-sided selfadjoint ideal. Let $\mathsf{B}:= \mathsf{A}/N_{\phi}$ and let $\pi_{\phi}:\mathsf{A} \to \mathsf{B}$ denote
the canonical quotient map. We want to define the coproduct on $\mathsf{B}$ by the formula
\begin{equation} \label{copb} \Delta_{\mathsf{B}}
\circ \pi_{\phi}(a) = (\pi_{\phi} \otimes \pi_{\phi})\circ \Delta(a) \in \mathsf{B} \otimes \mathsf{B}.
\end{equation}
We need to check that it is well-defined - to this end we employ a slightly modified idea from the proof of Theorem 2.1 of
\cite{bedos+murphy+tuset01}. A standard use of Cauchy-Schwarz inequality implies that $\phi|_{N_{\phi}}=0$, so that there exists
a faithful state $\psi \in \mathsf{B}^*$ such that $\psi \circ \pi_{\phi} = \phi$. Faithfulness of $\psi$ implies that also the
map ${\rm id}_{\mathsf{B}} \otimes \psi:\mathsf{B} \otimes \mathsf{B} \to \mathsf{B}$ is faithful (note that here faithfulness of
a positive map $T$ is understood in the usual sense, namely $Ta = 0$ and $a\geq 0$ imply $a=0$) and thus also $\psi \otimes
\psi\in (\mathsf{B} \otimes \mathsf{B})^*$ is faithful. Suppose then that $a \in N_{\phi}$. We have then
\[
0 = \phi(a^*a) =
(\phi \otimes \phi)\circ\Delta(a^*a) = (\psi \otimes \psi)\circ (\pi_{\phi} \otimes \pi_{\phi})(\Delta(a^*a)),
\]
so also $(\pi_{\phi} \otimes \pi_{\phi})
\Delta(a^*a) = 0$. The last statement implies that $(\pi_{\phi} \otimes \pi_{\phi}) \Delta(a) = 0$ and validity of the definition given
in the formula \eqref{copb} is established. The fact that $\Delta_{\mathsf{B}}$ is a coassociative unital $*$-homomorphism follows
immediately from the analogous properties of $\Delta$; similarly the cancellation properties of $\mathsf{B}$ follow from obvious
equalities of the type
\[
(\mathsf{B} \otimes 1_{\mathsf{B}})\Delta_{\mathsf{B}}(\mathsf{B})  =
(\pi_{\phi} \otimes \pi_{\phi}) ((\mathsf{A} \otimes
1_{\mathsf{A}})\Delta(\mathsf{A}))
\]
and the cancellation properties of $\mathsf{A}$. Thus $(\mathsf{B}, \Delta_{\mathsf{B}})$ is the algebra of continuous functions on a compact quantum group $\QG'$ and it remains to check that
$\psi$ defined above is actually the invariant state on $\mathsf{B}$. This is however an immediate consequence of the following
observation:
\[
(\psi \otimes \psi)\circ\Delta_{\mathsf{B}}\circ\pi_{\phi} =(\phi \otimes \phi
) \circ \Delta = \phi = \psi\circ \pi_{\phi},
\]
so that $\psi$ is an idempotent state and, as it is faithful, it has to coincide with the Haar state of $\QG'$ (\cite{woronowicz98}).
\end{proof}

Note that in fact the theorem remains valid without the assumption of faithfulness of $h_{\QG'}$. The proof of the `only
if' part remains the same, and the other implication follows from the modular properties of Haar states on
not-necessarily-coamenable compact quantum groups implying that their null spaces are always selfadjoint.

The following proposition proved in  \cite[Section 3]{franz+skalski08a} will be useful for the classification of idempotent
states in the next two sections.

\begin{proposition}\label{prop-antipode}
Let $\phi\in C(\QG)^*$ be an idempotent state. Then $\phi$ is invariant
under the antipode, in the sense that $\phi(a)=\phi\circ S(a)$ for all $a\in \Pol$.
\end{proposition}

Note that the states invariant under the antipode are automatically invariant under the \emph{scaling automorphism group}
$\{\tau_t:t \in \mathbb{R}\}$ discussed in \cite{woronowicz98}. This is the content of the next proposition:

\begin{proposition} Let $\omega \in C(\QG)^*$ satisfy the condition $\omega|_{\Pol} = \omega \circ
S|_{\Pol}$. Then
\[\omega \circ \tau_t = \omega, \;\;\; t \in \mathbb{R}.\]
\end{proposition}

\begin{proof}
The assumed invariance of $\omega$ under the antipode is equivalent to the equality $\omega \circ \tau_{-i}|_{\Pol} =
\omega|_{\Pol}$. The idea of the rest of the proof is based on the use of Woronowicz characters, as in the proof of Lemma
2.9(2) in \cite{tomatsu07}. Let $v \in B(\mathsf{H}_v) \odot \Pol $ be an irreducible unitary representation of $\QG$.
We have
\[ T_v:=(\textup{id}_{B(\mathsf{H}_v)} \otimes \omega)(v) = (\textup{id}_{B(\mathsf{H}_v)} \otimes \omega \circ \tau_{-i})(v) =
F_v (\textup{id}_{B(\mathsf{H}_v)} \otimes \omega)(v)  (F_v)^{-1},\] where $F_v= (\textup{id}_{B(\mathsf{H}_v)} 
\otimes f_1)(v)$ and $f_1$ is the Woronowicz character on $\Pol$. The last formula implies that $T_v$ commutes 
with $F_v$. As $F_v$ is a (strictly) positive operator on a finite dimensional Hilbert space, we can consider the 
unitary group it generates and deduce immediately that for each $t \in \mathbb{R}$ we have $F_v^{it} T_v = T_v 
F_v^{it}$. But this implies that
\[ (\textup{id}_{B(\mathsf{H}_v)} \otimes \omega)(v) = F_v^{it} (\textup{id}_{B(\mathsf{H}_v)} \otimes \omega)(v)  F_v^{-it} =
(\textup{id}_{B(\mathsf{H}_v)} \otimes \omega\circ \tau_t)(v).\] The fact that $\Pol$ is a linear span of the coefficients
of the irreducible unitary representations of $\QG$ shows that $\omega(a) = \omega (\tau_t(a))$ for each $a\in \Pol$;
density of $\Pol$ in $C(\QG)$ ends the proof.
\end{proof}

\begin{corollary}
Let $\phi \in C(\QG)^*$ be an idempotent state. Then $\phi$ is preserved by the scaling automorphism group and 
the conditional expectation $E_{\phi}:=(\phi \otimes \textup{id}_{C(\QG)}) \circ \Delta$ associated to  $\phi$ 
commutes with the modular group of the Haar state:
\[ E_{\phi} \circ \sigma_t = \sigma_t \circ E_{\phi}, \;\;\; t \in \mathbb{R}.\]
\end{corollary}

\begin{proof}
It is an immediate consequence of the last lemma and the commutation relation $\Delta \circ \sigma_t = (\sigma_t \otimes
\tau_{-t}) \circ \Delta$ ($t \in \mathbb{R}$). The fact that $E_{\phi}$ is a conditional expectation was established in
\cite{franz+skalski09}.
\end{proof}

\subsection{Idempotent states on cocommutative compact quantum groups}

Suppose now that a compact quantum group $\QG$ is cocommutative, i.e.\ $\Delta = \tau\circ \Delta$, where $\tau: 
C(\QG) \otimes C(\QG) \to C(\QG) \otimes C(\QG)$ denoted the usual tensor flip. It is easy to deduce from the 
general theory of duality for quantum groups (\cite{kustermans+vaes00}) that $C(\QG)$ is isomorphic to the 
$C^*$-algebra of a (classical) discrete group $\Gamma$, which should be thought of as the algebra of continuous 
functions on a quantum group dual $\hat{\Gamma}$. Note that the notation $C^*(\Gamma) \approx C(\hat{\Gamma})$, 
which can be considered as a definition of the compact quantum group $\hat{\Gamma}$ is compatible with the usual 
Pontriagin duality for locally compact abelian groups. For the reasons mentioned in the previous section in 
general we need to distinguish between the reduced and the universal version of $C(\hat{\Gamma})$; thus we 
restrict our attention to amenable $\Gamma$. The following generalises Theorem 6.2 of \cite{franz+skalski08a} to 
the infinite-dimensional context.

\begin{theorem} \label{cocom}
Let $\Gamma$ be an amenable discrete group and $\mathsf{A}=C^*(\Gamma)\approx C(\hat{\Gamma})$. There is a one-to-one correspondence between idempotent states
on $\mathsf{A}$ and subgroups of $\Gamma$. An idempotent state $\phi\in \mathsf{A}^*$ is a Haar idempotent if and only if the corresponding
subgroup of $\Gamma$ is normal.
\end{theorem}
\begin{proof}
The dual of $\mathsf{A}$ may be identified with the Fourier-Stieltjes algebra $B(\Gamma)$. The convolution of functionals in $\mathsf{A}^*$
corresponds then to the pointwise multiplication of functions in $B(\Gamma)$ and $\phi\in B(\Gamma)$ corresponds to a positive
(respectively, unital) functional on $\mathsf{A}$ if and only if it is positive definite (respectively, $\phi(e)=1$). This implies that
$\phi \in B(\Gamma)$ corresponds to an idempotent state if and only if it is an indicator function (of a certain subset $S
\subset \Gamma$) which is positive definite. It is a well known fact that this happens if and only if $S$ is a subgroup of
$\Gamma$ (\cite{hewitt+ross70}, Cor. (32.7) and Example (34.3 a)). It remains to prove that if $S$ is a subgroup of $\Gamma$ then
$\chi_S \in B(\Gamma)$ is a Haar state on a quantum subgroup of $\hat{\Gamma}$ if and only if $S$ is normal. For the `if' direction
assume that $S$ is a normal subgroup and consider the $C^*$-bialgebra $\mathsf{B}=C^*(\Gamma/S)$ (recall that quotients of
amenable groups are amenable). Let $F(\Gamma)$ denote the dense $*$-subalgebra of $\mathsf{A}$ given by the functions $f= \sum_{\gamma
\in \Gamma} \alpha_{\gamma} \lambda_{\gamma}$ ($\alpha_{\gamma} \in \mathbb{C}$, $\{\gamma: \alpha_{\gamma} \neq 0\}$ finite). Define
$j:F(\Gamma) \to \mathsf{B}$ by
\[
j(f)= \sum_{\gamma \in \Gamma}  \alpha_{\gamma}\lambda_{[\gamma]},
\]
where $f$ is as above. So-defined $j$ is bounded: note that it is a restriction of the transpose of the map
$T: \mathsf{B}^* \to \mathsf{A}^*$ given by
\[
T(\phi) (\gamma) = \phi( [\gamma]), \; \phi \in B(\Gamma / S),  \gamma \in \Gamma.
\]
The map $T$ is well defined as it maps positive definite functions into positive definite functions; these generate the relevant
Fourier-Stieltjes algebras. Further the closed graph theorem allows to prove that $T$ is bounded; therefore so is $T^*: \mathsf{A}^{**}
\to \mathsf{B}^{**}$ and $j=T^*|_{F(\Gamma)}$. It is now easy to check that the extension of $j$ to $\mathsf{A}$ is a surjective unital
$*$-homomorphism (onto $\mathsf{B}$).
 As the invariant state on $\mathsf{B}$ is given by
\[
h_{\mathsf{B}} \left(\sum_{\kappa \in \Gamma/S}  \alpha_{\kappa} \lambda_{\kappa} \right) = \alpha_{[e]},
\]
there is
\[
h_{\mathsf{B}} (j(f)) = \sum_{\gamma \in S} \alpha_{\gamma},
\]
so that $h_{\mathsf{B}} \circ j$ corresponds via the identification of $\mathsf{A}^*$ to $B(\Gamma)$
exactly to the characteristic function of $S$.

The other direction follows exactly as in \cite{franz+skalski08a}; we reproduce the argument for the sake of completeness.
    Suppose that $S$ is a subgroup of $\Gamma$ which is not normal and let $\gamma_0\in \Gamma$, $s_0 \in S$ be such that $\gamma_0 s_0 \gamma_0 ^{-1} \notin S $. Denote by $\phi_S$ the state on $\mathsf{A}$ corresponding to the indicator function of $S$. Define $f \in \mathsf{A}$ by
$f= \lambda_{\gamma_0 s_0} - \lambda_{\gamma_0}$. Then
\[
f^*f = 2 \lambda_e - \lambda_{s_0^{-1}} - \lambda_{s_0}, \;\;\;\; f f^* =  2 \lambda_e - \lambda_{\gamma_0 s_0^{-1}\gamma_0^{-1}} - \lambda_{\gamma_0 s_0 \gamma_0^{-1}}.
\]
This implies that
\[
\phi_S(f^*f) = 0,\;\;\; \phi_S(f f^*) = 2,
\]
so that ${\rm Ker}\, \phi_S$ is not selfadjoint and $\phi_S$ cannot be a Haar idempotent.
\end{proof}

\begin{corollary}
Let $\QG$ be a  coamenable cocommutative compact quantum group. The following are equivalent:
\begin{enumerate}
\item all idempotent states on $\QG$ are Haar idempotents; \item $C(\QG)\cong C^*(\Gamma)$ for an amenable 
Dedekind (i.e.\ containing no non-normal subgroups) discrete group $\Gamma$.
\end{enumerate}
\end{corollary}

Finite Dedekind groups have been first studied by Dedekind and then characterised by Baer in \cite{baer66}. They 
are also sometimes called \emph{quasi-Hamiltonian} (and \emph{Hamiltonian} groups are non-abelian Dedekind 
groups).

\begin{remark}
If $\Gamma$ is a discrete abelian group, then there is a one-to-one correspondence between subgroups of 
$\hat{\Gamma}$ (the Pontriagin dual of $\Gamma$) and subgroups of $\Gamma$: if $N$ is a subgroup of $\Gamma$, 
then $\widehat{\Gamma/N}$ is a subgroup of $\hat{\Gamma}$. This is a part of a general fact that the subgroup of 
a dual group corresponds to the dual of a quotient group. Theorem \ref{cocom} and its proof show that the 
analogous statement remains true for not necessary abelian discrete groups $\Gamma$ -- quantum subgroups of 
$\hat{\Gamma}$ are given by $\widehat{\Gamma/N}$, where $N$ is a \emph{normal} subgroup of $\Gamma$ (normality is 
needed to make $\Gamma/N$ a group). It is natural to ask whether one can construct  the $C^*$-algebra 
$C^*(\Gamma/S)$ and equip it with some extra algebraic structure if $S$ is a non-normal subgroup of $\Gamma$. A 
partial answer to this question (for finite $S$) can be found in \cite{DelvauxVanDaele}.
\end{remark}

\section{Idempotent states on $U_q(2)$ ($q\in(-1,0)\cup(0,1]$)}

For $q=1$, $C(U_q(2))$ is equal to the $C^*$-algebra of continuous functions
on the unitary group $U(2)$, and by Kawada and It\^o's classical theorem all
idempotent states on $C(U(2))$ come from Haar measures of compact subgroups of
$U(2)$. In this section we shall classify the idempotent states on $C(U_q(2))$
for $-1<q<1, q \neq 0$. It turns out that they all correspond to Haar
states of quantum subgroups of $U_q(2)$.

We begin with some preparatory lemmas.

\begin{lemma}\label{lem-diag}
Let $\phi:{\rm Pol}(U_q(2))\to\mathbb{C}$ be an idempotent state. Then we have
\[
\phi(u^{(s)}_{k\ell} v^r) = 0 \qquad\mbox{ if } \qquad k\not=\ell,
\]
and $\phi(u^{(s)}_{kk} v^r) \in \{0,1\}$,
for all $s\in\frac{1}{2}\mathbb{Z}_+$, $r\in \mathbb{Z}$, $-s\le k,\ell\le s$.
\end{lemma}
\begin{proof}
By Proposition \ref{prop-antipode}, we have $\phi\circ S=\phi$ on ${\rm
  Pol} (U_q(2))$. Therefore, by Equation \eqref{u2-s-squared},
\[
\phi(u^{(s)}_{k\ell} v^r) = \phi\circ S^2(u^{(s)}_{k\ell} v^r) = q^{2(k-\ell)} \phi(u^{(s)}_{k\ell} v^r)
\]
i.e.\ $\phi(u^{(s)}_{k\ell} v^r)=0$ for $k\not=\ell$.

Define the matrices $M_{s,p}(\phi)\in M_{2s+1}(\mathbb{C})$ by
\[
M_{s,p}(\phi) = \left(\phi(u^{(s)}_{k\ell} v^{p+s+\ell})\right)_{-s\le
  k,\ell\le s}.
\]
Then $\phi=\phi\star\phi$ is equivalent to
\[
M_{s,p}(\phi) = \big(M_{s,p}(\phi)\big)^2
\]
for all $s\in\frac{1}{2}\mathbb{Z}_+$, $p\in \mathbb{Z}$. As we have already
seen that these matrices are diagonal, it follows that the diagonal entries
can take only the values $0$ and $1$.
\end{proof}

\begin{lemma}\label{lem-torus}
If $\phi:{\rm Pol} (U_q(2))\to\mathbb{C}$ is an idempotent state with
$\phi(u^{(1)}_{00})=1$, then there exists an idempotent state
$\tilde{\phi}:{\rm Pol}(\mathbb{T}^2)\to\mathbb{C}$ such that
\[
\phi = \tilde{\phi}\circ \pi_{\mathbb{T}^2}.
\]
\end{lemma}
\begin{proof}
By the previous lemma $\phi(c)=\phi(c^*)=0$

We have  $u^{(1)}_{00}=1-(1+q^2) c^*c$, therefore
$\phi(u^{(s)}_{00})=1$ is equivalent to $\phi(c^*c)=0$. Then by
Theorem  \ref{mult-dom}, $c,c^*\in D_\phi$, and $\phi$ vanishes on expressions of the form $uc$, $c u$, $uc^*$, $c^* u$ with $u\in {\rm Pol}(U_q(2))$. But since $vc =c v$ and $u^{(s)}_{k\ell}c = q^{-(k-\ell)}c u^{(s)}_{k\ell}$, $u^{(s)}_{k\ell}c^* = q^{-(k-\ell)}c^* u^{(s)}_{k\ell}$ for $s\in\frac{1}{2}\mathbb{Z}_+$, $-s\le k,\ell\le s$, we can deduce that $\phi$ vanishes on the ideal
\[
\mathcal{I}_c = \{ u_1c u_2,u_1c^* u_2 ; u_1,u_2\in {\rm
  Pol}(U_q(2))\}
\]
generated by $c$ and $c^*$. It follows that we can divide
out $\mathcal{I}_c$, i.e.\ there exists a unique state $\tilde{\phi}$ on
${\rm Pol}(U_q(2))/ \mathcal{I}_c$ such that the diagram
\[
\xymatrix{
{\rm Pol}(U_q(2))\ar[r]^\pi \ar[d]_\phi & {\rm Pol}(U_q(2))/ \mathcal{I}_c \ar[dl]^{\tilde{\phi}} \\
\mathbb{C}
}
\]
commutes.

But $\varepsilon(\mathcal{I}_c)=0$,
\[
\Delta(\mathcal{I}_c) \subseteq \mathcal{I}_c\odot{\rm
  Pol}(U_q(2))+{\rm Pol}(U_q(2)) \odot \mathcal{I}_c,
\]
and $S(\mathcal{I}_c)\subseteq \mathcal{I}_c$
i.e.\ $\mathcal{I}_c$ is also a Hopf $*$-ideal and ${\rm Pol}(U_q(2))/
\mathcal{I}_c$ is a $*$-Hopf algebra. One easily verifies that actually
${\rm Pol}(U_q(2))/ \mathcal{I}_c\cong {\rm Pol}(\mathbb{T}^2)$. Since
$\pi_{\mathbb{T}^2}:{\rm Pol}(U_q(2))\to {\rm Pol}(\mathbb{T}^2)$ is
surjective coalgebra
morphism, its dual $\pi_{\mathbb{T}^2}^*:({\rm Pol}(\mathbb{T}^2))^*\ni
f\mapsto \pi_{\mathbb{T}^2}^*(f)=f\circ\pi_{\mathbb{T}^2}\in ({\rm
  Pol}(U_q(2)))^*$ is an injective algebra homomorphism, and
$\tilde{\phi}=(\pi_{\mathbb{T}^2}^*)^{-1}(\phi)$ is again idempotent.
\end{proof}

\begin{lemma}\label{lem-haar1}
If $\phi:{\rm Pol}(U_q(2))\to\mathbb{C}$ is an idempotent state with
$\phi(u^{(1)}_{00})=0$, then $\phi(u^{(s)}_{00})=0$ for all integers $s\ge 1$,
i.e.\ we have $\phi|_{\mathbb{C}[c^*c]}=h|_{\mathbb{C}[c^*c]}$.
\end{lemma}
\begin{proof}
Recall $u^{(1)}_{00}=1-(1+q^2) c^*c$. Therefore $\phi(u^{(s)}_{00})=0$ is equivalent to $\phi(c^*c)=\frac{1}{1+q^2}$.

Assume there exists an integer $s > 1$ with $\phi(u^{(s)}_{00})=1$. Then the Cauchy-Schwarz
inequality implies $\phi((u^{(s)}_{00})^*u^{(s)}_{00})\ge 1$. The unitarity of the representation
$v^{(s,p)}$ gives
\[
1=\sum_{k=-s}^s \left(u^{(s)}_{k0}v^{p+s}\right)^*u^{(s)}_{k0}v^{p+s} = \sum_{k=-s}^s \left(u^{(s)}_{k0}\right)^*u^{(s)}_{k0},
\]
therefore
\[
\phi\left(\sum_{\genfrac{.}{.}{0pt}{}{k\in\{-s,\ldots,s\}}{k\not=0}}\left(u^{(s)}_{k0}\right)^*u^{(s)}_{k0}\right)
\le 0,
\]
and in particular $\phi\big((u^{(s)}_{s0})^*u^{(s)}_{s0}\big)=0$.
We have
\begin{eqnarray*}
\left(u^{(s)}_{s}\right)^*u^{(s)}_{s0} &=& \left[\begin{array}{c} 2s \\ s
  \end{array}\right]_{q^2}\big((a^*)^sc^s\big)^*(a^*)^sc^s = \left[\begin{array}{c} 2s \\ s
  \end{array}\right]_{q^2}(c^*)^sc^sa^s(a^*)^s \\
&=&(c^*c)^s (1-q^2c^*c)\cdots (1-q^{2s}c^*c)
\end{eqnarray*}
By the representation theory of $C(SU_q(2))$, $c^*c$ is positive self-adjoint
contraction, with the spectrum $\sigma(c^*c)\subseteq\{q^{2n};n\in\mathbb{Z}_+\}\cup\{0\}$, and therefore the
product $(1-q^2c^*c)\cdots (1-q^{2s}c^*c)$ defines a strictly
positive operator. Therefore
$\phi\left((u^{(s)}_{s,0})^*u^{(s)}_{s,0}\right) = 0$ implies
$\phi\big((c^*c)^s\big)=0$, which is impossible if $\phi(c^*c)=\frac{1}{1+q^2}>0$.

Therefore $\phi(u^{(s)}_{00})=0$ for all
integers $s\ge 1$.
\end{proof}

\begin{lemma}\label{lem-haar2}
Let
\[
\mathcal{A}_0 = {\rm span}\{u^{(s)}_{k\ell}v^r;s\in\frac{1}{2}\mathbb{Z}_+,
s>0, -s\le k,\ell\le s,r\in\mathbb{Z}\}
\]
i.e.\ $\mathcal{A}_0$ is the subspace spanned by the matrix coefficients of the unitary irreducible 
representations of dimension at least two.

Assume that $\phi|_{\mathbb{C}[c^*c]}=
h|_{\mathbb{C}[c^*c]}$, i.e.\
\[
\phi(u^{(s)}_{00})=\delta_{0s}
\]
for $s\in\mathbb{Z}_+$.

Then we have $\phi|_{\mathcal{A}_0}=h|_{\mathcal{A}_0} $, i.e.
\[
\phi\left(u^{(s)}_{k\ell}v^r\right)=0
\]
for all $r\in \mathbb{Z}$, $s\in\frac{1}{2}\mathbb{Z}_+$, $s>0$, and $-s\le
k,\ell\le s$.
\end{lemma}
\begin{proof}
By Lemma \ref{lem-diag}, we already know that
$\phi(u^{(s)}_{k\ell}v^r)\in\{0,1\}$, and $\phi(u^{(s)}_{k\ell}v^r)=0$  for
$k\not=\ell$. Assume there exist $s\in\frac{1}{2}\mathbb{Z}_+$, $s>0$, $-s\le
k\le s$ and $r\in\mathbb{Z}$ such that
\[
\phi(u^{(s)}_{kk}v^r)=1.
\]
We will show that this is impossible, if $\phi$ agrees with the Haar state $h$
on the subalgebra generated by $c^*c$.

By the Cauchy-Schwarz inequality, we have
\[
\phi\left((u^{(s)}_{kk})^*u^{(s)}_{kk}\right)\ge
\left|\phi(u^{(s)}_{kk}v^r)\right|^2 = 1.
\]
Applying $\phi$ to
\[
\sum_{\ell=-s}^s (u^{(s)}_{\ell k})^*u^{(s)}_{\ell k} =1,
\]
we can deduce
\[
\sum_{\genfrac{.}{.}{0pt}{}{\ell\in\{-s,\ldots,s\}}{\ell\not=k}} \phi\left((u^{(s)}_{\ell
  k})^*u^{(s)}_{\ell k}\right) =0.
\]
But this contradicts $\phi|_{\mathbb{C}[c^*c]}=
h|_{\mathbb{C}[c^*c]}$, because
\[
\sum_{\genfrac{.}{.}{0pt}{}{\ell\in\{-s,\ldots,s\}}{\ell\not=k}}
(u^{(s)}_{\ell
  k})^*u^{(s)}_{\ell k} = 1- (u^{(s)}_{kk})^*u^{(s)}_{kk}
\]
is a non-zero positive element in $\mathbb{C}[c^*c]$ and the Haar
state is faithful.
\end{proof}

Before we formulate the main theorem of this section we need two more remarks which will be used in  the proof.

\begin{remark}
As $C(\mathbb{T}^2)$ is commutative, by Kawada and
It\^o's theorem all idempotent states on $C(\mathbb{T}^2)$ are induced by Haar
measures of compact subgroups of the two-dimensional torus $\mathbb{T}^2$.
\end{remark}

\begin{remark}
As a compact quantum group, $\mathbb{Z}_n$ is given by
\[
{\rm Pol}(\mathbb{Z}_n)=C(\mathbb{Z}_n)={\rm
  span}\,\{w_0,\ldots,w_{n-1}\},
\]
with $w_kw_\ell =
w_{k+\ell\!\mod n}$, $S(w_k)=w_{n-k}=(w_k)^*$, $\Delta(w_k)=w_k\otimes w_k$, and
$\varepsilon(w_k)=1$ for $k=0,\ldots, n-1$. The Haar state of $\mathbb{Z}_n$
is given by $h(w_k)=\delta_{0k}$. $C(\mathbb{Z}_n)$ can also be obtained from ${\rm
  Pol}(U(1))$ by dividing out the Hopf $*$-ideal $\{u_1(w^n-1)u_2;\,
u_1,u_2\in {\rm Pol}(U(1))\}$.

Analogous to \cite[Section 4]{wysoczanski04},
one can define the twisted product $SU_q(2)\ltimes_\sigma
\mathbb{Z}_n$. Alternatively, ${\rm Pol}(SU_q(2)\ltimes_\sigma \mathbb{Z}_n)$
can be obtained from $U_q(2)\cong SU_q(2)\ltimes_\sigma U(1)$ by
dividing out the Hopf ideal $\{u_1(v^n-1)u_2;\,u_1,u_2\in {\rm
  Pol}(U_q(2))\}$, and $C(SU_q(2)\ltimes_\sigma \mathbb{Z}_n)$ as its
$C^*$-completion. This construction shows that $SU_q(2)\ltimes_\sigma \mathbb{Z}_n$ is a quantum subgroup of 
$U_q(2)\cong SU_q(2)\ltimes_\sigma U(1)$. As in the case of $U_q(2)$, $C(SU_q(2)\ltimes_\sigma 
\mathbb{Z}_n)=C(SU_q(2))\otimes C(\mathbb{Z}_n)$ and the Haar state of $SU_q(2)\ltimes_\sigma \mathbb{Z}_n$ is 
equal to the tensor product of the Haar states of $SU_q(2)$ and $\mathbb{Z}_n$.
\end{remark}

We can now give a description of all idempotent states on $U_q(2)$. It turns
out that they are all induced by Haar states of quantum subgroups of $U_q(2)$.

\begin{theorem}\label{theo-uq2}
Let $q\in(-1,0)\cup(0,1)$. Then the following is a complete list of the idempotent states on the
compact quantum group $U_q(2)$.
\begin{enumerate}
\item
The Haar state $h$ of $U_q(2)$.
\item
$\tilde{\phi}\circ\pi_{\mathbb{T}^2}$, where $\pi_{\mathbb{T}^2}$ denotes the
  surjective quantum group morphism $\pi_{\mathbb{T}^2}:C(U_q(2))\to
    C(\mathbb{T}^2)$ and $\tilde{\phi}$ is an idempotent state on
    $C(\mathbb{T}^2)$. In particular if $\widetilde{\varepsilon}$ denotes the counit of  $C(\mathbb{T}^2)$, then
   $\widetilde{\varepsilon}\circ\pi_{\mathbb{T}^2}$ is the counit of $U_q(2)$.
\item
The states induced by the Haar states of the compact quantum subgroups $SU_q(2)\ltimes_\sigma
\mathbb{Z}_n$ of $U_q(2)\cong SU_q(2)\ltimes_\sigma U(1)$, for $n\in \mathbb{N}$. The case of $n
=1$ corresponds to the Haar state on $SU_q(2)$ viewed as a quantum subgroup of $U_q(2)$.
\end{enumerate}
\end{theorem}

\begin{proof} 
Let $\phi:C(U_q(2))\to \mathbb{C}$ be an idempotent state on $U_q(2)$. Clearly
$\phi$ is uniquely determined by its restriction to ${\rm Pol}(U_q(2))$.

We distinguish two cases.

\paragraph{Case (i) $\phi(u^{(1)}_{00})=1$}
In this case Lemma \ref{lem-torus} shows that $\phi$ is induced by an idempotent
state on the quantum subgroup $\mathbb{T}^2$ of $U_q(2)$, i.e.\
$\phi=\tilde{\phi}\circ\pi_{\mathbb{T}^2}$ for some idempotent state
$\tilde{\phi}:C(\mathbb{T}^2)\to\mathbb{C}$. This case includes the counit
$\varepsilon$ of $U_q(2)$, it corresponds to the trivial subgroup $\{1\}$ of
$\mathbb{T}^2$.

\paragraph{Case (ii) $\phi(u^{(1)}_{00})=0$}
In this case Lemma \ref{lem-haar1} and Lemma \ref{lem-haar2} imply that $\phi$
agrees with the Haar state $h$ on the subspace $\mathcal{A}_0$, i.e.\
\[
\phi(u^{(s)}_{k\ell} v^r) =0
\]
for all $s\in\frac{1}{2}\mathbb{Z}_+$, $s>0$, $-s\le k,\ell\le s$, and
$r\in\mathbb{Z}$. It remains to determine $\phi$ on the $*$-subalgebra
 ${\rm alg}\{v,v^*\}$ generated by $v$, since ${\rm
  Pol}(U_q(2))=\mathcal{A}_0\oplus {\rm alg}\{v,v^*\}$ as a vector space. But this subalgebra is isomorphic to
the $*$-Hopf algebra ${\rm Pol}(U(1))$ of
polynomials on the unit circle, and therefore $\phi|_{{\rm alg}\{v,v^*\}}$ has
to be induced by the Haar measure of a compact subgroup of $U(1)$. We have the
following possibilities.
\begin{enumerate}
\item
$\phi|_{{\rm alg}\{v,v^*\}}=\varepsilon_{U(1)}$, i.e.\ the restriction of
$\phi$ to ${\rm alg}\{v,v^*\}$ is equal to the counit of ${\rm Pol}(U(1))$ . In this case we have
\[
\phi(u^{(s)}_{k\ell} v^r) =\left\{
\begin{array}{cll}
1 & \mbox{ if } & s=k=\ell=0,\, \mbox{ and } r\in\mathbb{Z}, \\
0 & \mbox { else}. &
\end{array}\right.
\]
This formula shows that $\phi=h_{SU_q(2)}\circ \pi_{SU_q(2)}$, where
$\pi_{SU_q(2)}$ is the quantum groups morphism from $C(U_q(2))$ onto $C(SU_q(2))$
and $h_{SU_q(2)}$ denotes the Haar state of $SU_q(2)$.
\item
$\phi|_{{\rm alg}\{v,v^*\}}=h_{U(1)}$, i.e.\ the restriction of
$\phi$ to ${\rm alg}\{v,v^*\}$ is equal to the Haar state of ${\rm Pol}(U(1))$ . In this case we have
\[
\phi(u^{(s)}_{k\ell} v^r) =\left\{
\begin{array}{cll}
1 & \mbox{ if } & s=k=\ell=0 \mbox{ and } r=0, \\
0 & \mbox { else}. &
\end{array}\right.
\]
We see that in this case $\phi$ is the Haar state $h$ of $U_q(2)$.
\item
$\phi|_{{\rm alg}\{v,v^*\}}$ is the idempotent state on $U(1)$ induced by the
Haar measure of the subgroup $\mathbb{Z}_n\subseteq U(1)$ for some $n\in\mathbb{N}$,
$n\ge 2$. In this case we have
\[
\phi(u^{(s)}_{k\ell} v^r) =\left\{
\begin{array}{cll}
1 & \mbox{ if } & s=k=\ell=0 \mbox{ and } r\equiv 0\mod n, \\
0 & \mbox { else}. &
\end{array}\right.
\]
It follows that $\phi$ is induced by the Haar state of the quantum subgroup
$SU_q(2)\ltimes_\sigma \mathbb{Z}_n$ of $U_q(2)\cong SU_q(2)\ltimes_\sigma U(1)$.
\end{enumerate}
Conversely, all the states we have found are induced by Haar states on
quantum subgroups of $U_q(2)$, therefore they are clearly idempotent. It can be also checked directly.
\end{proof}

We see that all idempotent states on $U_q(2)$ are induced from Haar states of quantum subgroups of $U_q(2)$. We can also deduce the complete list of quantum subgroups of $U_q(2)$.

\begin{corollary}\label{uq2-subgroups}
Let $q\in(-1,0)\cup(0,1)$. Then the following is a complete list of the non-trivial quantum subgroups of $U_q(2)$.
\begin{enumerate}
\item
The two-dimensional torus and its closed subgroups.
\item
The compact quantum groups of the form $SU_q(2)\ltimes_\sigma \mathbb{Z}_n$, with $n\in\mathbb{N}$
(here the twisting is identical to that appearing in the identification $U_q(2) \cong
SU_q(2)\ltimes_\sigma \mathbb{T}$).
\end{enumerate}
\end{corollary}

\section{Idempotent states on compact quantum groups $SU_q(2)$ and $SO_q(3)$ ($q \in (-1,0)\cup(0,1]$)}

Let us first discuss the case $q=1$. $C(SU_1(2))$ and $C(SO_1(3))$ are the
algebras of continuous functions on the groups $SU(2)$ and $SO(3)$. All
idempotent states correspond  to Haar measures on compact subgroups. The list of
these subgroups can be found, e.g., in \cite{podles95}.

Consider now the generic case $q\in(-1,0)\cup(0,1)$.
Every idempotent state on $SU_q(2)$ induces an idempotent state on $U_q(2)$,
since $SU_q(2)$ is a quantum subgroup of $U_q(2)$. This observation allows us to deduce all idempotent states on $SU_q(2)$ from Theorem \ref{theo-uq2}. We omit
the details and just state the result.

\begin{theorem}\label{theo-suq2}
Let $q\in(-1,0)\cup(0,1)$. The Haar state, the counit, and the idempotent states induced by the quantum
subgroups $U(1)$ and $\mathbb{Z}_n$, $2\le n\le \infty$, are the only
idempotent states on $SU_q(2)$.
\end{theorem}

Since the morphism $j:C(SU_q(2))\to C(U(1))$ gives the diagonal matrices
\[
\left(j(u^{(s)}_{k\ell})\right)_{-s\le k,\ell\le s} =
\left(\begin{array}{cccc}
z^{2s} &  & & \\
 & z^{2s-2} & & \\
 & & \ddots & \\
 & & & z^{-2s}
\end{array}\right),
\]
we get
\begin{equation} \label{u1-idemp}
\big(h_{U(1)}\circ j\big)(u^{(s)}_{k\ell}) = \left\{
\begin{array}{lll}
1 & \mbox{ if } & s\in\mathbb{Z}_+,\, k=\ell=0, \\
0 & \mbox{ else.}&
\end{array}\right.
\end{equation}
and
\begin{eqnarray}\label{z2n-idemp}
\big(h_{\mathbb{Z}_{2n}}\circ j\big)(u^{(s)}_{k\ell}) &=& \left\{
\begin{array}{lll}
1 & \mbox{ if } & s\in\mathbb{Z}_+,\,\, k=\ell,\, 2k\equiv0 \!\!\mod 2n, \\
0 & \mbox{ else.}&
\end{array}\right. \\
\label{z2n+1-idemp}
\big(h_{\mathbb{Z}_{2n+1}}\circ j\big)(u^{(s)}_{k\ell}) &=& \left\{
\begin{array}{lll}
1 & \mbox{ if } & k=\ell, \,\, 2k\equiv0 \!\!\mod 2n+1, \\
0 & \mbox{ else.}&
\end{array}\right.
\end{eqnarray}
for $n\in\mathbb{N}$.

Consider now the idempotent states on $SO_q(3)$. Since $C(SO_q(3))$ is a
subalgebra of $C(SU_q(2))$ and since the inclusion map is a quantum group
morphism, every idempotent state on $SU_q(2)$ gives an idempotent state on
$SO_q(3)$ by restriction. We will show that all idempotent states on $SO_q(3)$
arise in this way. It follows that all idempotent states on $SO_q(3)$ are
induced from Haar states of quantum subgroups.
\begin{theorem}\label{theo-soq3}
Let $q\in(-1,0)\cup(0,1)$ and $n$ an odd integer. Then the  restrictions to $C(SO_q(3))$ of the idempotent
states $h_{\mathbb{Z}_{n}}\circ j$ and $h_{\mathbb{Z}_{2n}}\circ j$ coincide.

Furthermore, the Haar state, the counit, and the states induced from the Haar
states on the quantum subgroups $U(1)\cong SO(2)$ and its closed subgroups are the
only idempotent states on $SO_q(3)$.
\end{theorem}
\begin{proof}
The first statement follows from Equations \eqref{z2n-idemp} and \eqref{z2n+1-idemp}.

Let now $\phi$ be an idempotent state on $SO_q(3)$. Denote by $E$ the conditional expectation from
$C(SU_q(2))$ onto $C(SO_q(3))$ introduced in Paragraph \ref{subsub-cond-exp}. Then
$\hat{\phi}=\phi\circ E$ defines an idempotent state on $SU_q(2)$ such that
$\phi=\hat{\phi}|_{C(SO_q(3))}$. With this observation Theorem \ref{theo-soq3} follows immediately
from Theorem \ref{theo-suq2}.
\end{proof}

\begin{remark}
This method applies to quotient quantum groups in general. If $(\QG', j)$ is a normal quantum subgroup of $\QG$ (see Definition \ref{defnormal}) then all idempotent states on $C(\QG/\QG')$ arise as restrictions of
idempotent states on $C(\QG)$.
\end{remark}

As a corollary to Theorems \ref{theo-suq2} and  \ref{theo-soq3}, we recover Podle\'s' classification \cite{podles95} of the quantum subgroups of $SU_q(2)$ and $SO_q(3)$.
\begin{corollary}
Let $q\in(-1,0)\cup(0,1)$. Then $U(1)\cong SO(2)$ and its closed subgroups are the
only non-trivial quantum subgroups of both $SU_q(2)$ and $SO_q(3)$.
\end{corollary}

\section{Idempotent states on compact quantum semigroups $U_0(2)$, $SU_0(2)$ and $SO_0(3)$}

In this section we compute all idempotent states on $U_0(2)$, $SU_0(2)$ and $SO_0(3)$.
As in the cases $q\neq 0$ considered earlier we begin with the  $C^*$-bialgebra $C(U_0(2))$. Again we first need  some preparatory observations and lemmas.

Note that $\pi_{\mathbb{T}^2}:C(U_q(2))\to C(\mathbb{T}^2)$ is a well-defined
$*$-algebra and coalgebra morphism also for $q=0$, so $\mathbb{T}^2$ and its
compact subgroups induce idempotent states on $U_0(2)$.

For $q=0$ the algebraic relations of $a$ and $c$ become
\[
cc^*=c^*c,\qquad aa^*=1,\qquad ac=ac^*=0,\qquad a^*a=1-c^*c.
\]
As $a$ is a coisometry, we have a decreasing family of orthogonal projections $(a^*)^na^n$,
$n\in\mathbb{N}$ , which are group-like, i.e.\ $\Delta ((a^*)^na^n) = (a^*)^na^n\otimes
(a^*)^na^n$, and $c^*c=1-a^*a$ is also an orthogonal projection.

Denote by $M$ the unital semigroup $U(1)\times(\mathbb{Z}_+\cup\{\infty\})$
with the operation
\[
(z_1,n_1)\cdot (z_2,n_2) = (z_1z_2,\min(n_1,n_2))
\]
for $z_1,z_2\in U(1)$, $n_1,n_2\in \mathbb{Z}_+\cup\{\infty\}$. This is an abelian semigroup with unit element 
$e_M=(1,\infty)$. Equip $\mathbb{Z}_+\cup\{\infty\}$ with the topology in which a subset of 
$\mathbb{Z}_+\cup\{\infty\}$ is open if and only if it is either an arbitrary subset of $\mathbb{Z}_+$ or the 
complement of a finite subset of $\mathbb{Z}_+$ (i.e.\ $\mathbb{Z}_+\cup\{\infty\}$ is the 
one-point-compactification of $\mathbb{Z}_+$), and equip $M=U(1)\times(\mathbb{Z}_+\cup\{\infty\})$ with the 
product topology.

The $C^*$-bialgebra $C(M)$ will play an important role in this section.

\begin{lemma}\label{id-M}
A probability measure $\mu$ on $M$ is idempotent if and only if
there exists an $n\in\mathbb{Z}_+\cup\{\infty\}$ and
an idempotent probability $\rho$ on $U(1)$ such that $\mu=\rho\otimes \delta_{n}$.
\end{lemma}
\begin{proof}
Any probability on $M$ can be expressed as a sum $\mu=\sum_{n=0}^\infty
\rho_n\otimes \delta_n + \rho_\infty\otimes \delta_{\infty}$, where $\rho_n$,
$n\in\mathbb{Z}_+\cup\{\infty\}$ are uniquely determined positive measures on $U(1)$ with total mass
$\sum_{n=0}^\infty\rho_n(U(1))+\rho_\infty(U(1))=1$, and $\delta_n$
denotes the Dirac measure on $\mathbb{Z}_+\cup\{\infty\}$, i.e.\
\[
\delta_n(Q) = \left\{
\begin{array}{cll}
1 & \mbox{ if } & n \in Q, \\
0 & \mbox{ if } & n \not\in Q,
\end{array}\right.
\]
for $Q\subseteq \mathbb{Z}_+\cup\{\infty\}$. We have
\begin{eqnarray*}
(\delta_n\star\delta_m)(Q) &=&
(\delta_n\otimes\delta_m)\Big(\big\{(k,\ell)\in(\mathbb{Z}_+\cup\{\infty\})^2;\,\min(k,\ell)\in
Q\big\}\Big) \\
&=& \left\{\begin{array}{cll}
1 & \mbox{ if } & \min(n,m) \in Q, \\
0 & \mbox{ if } & \min(n,m) \not\in Q,
\end{array}\right.
\end{eqnarray*}
i.e.\ $\delta_n\star\delta_m=\delta_{\min(n,m)}$. Therefore
\[
\mu^{\star 2}= \rho_\infty^{\star 2}\otimes\delta_{\infty} +
\sum_{n=0}^\infty\left(\rho_n\star\left( \rho_n + 2\sum_{m=n+1}^\infty
  \rho_m + 2 \,\rho_\infty\right)\right)\otimes\delta_n.
\]
Clearly, if $\rho_\infty$ is an idempotent probability on $U(1)$ and
$\rho_n=0$ for $n\in\mathbb{Z}_+$,
then $\mu=\rho_\infty\otimes\delta_\infty$ is idempotent.

Assume now that
$\rho_\infty(U(1))<1$. Then there exists a unique $n\in\mathbb{Z}_+$ such that
$\sum_{m=n+1}^\infty \rho_m(U(1))<1$, $\sum_{m=n}^\infty
\rho_m(U(1))=1$ (i.e.\ $n$ is the biggest integer $m$ for which $\rho_{\ge
  m}=\rho_\infty+ \sum_{k=m}^\infty\rho_k$ is a probability). Let
$p=\rho_n(U(1))$. If $\mu$ is idempotent, then we have
\[
p = \mu(U(1)\times\{n\}) = \mu^{\star 2} (U(1)\times\{n\}) =\left(
\rho_n\star(2\rho_{\ge n}-\rho_n)\right)(U(1))=2p-p^2.
\]
Since $p=\rho_n(U(1))>0$ by the choice of $n$, we get $p=1$, i.e.\ $\rho_m=0$
for $m\not=n$ and $\rho_n$ is a probability. Then $\rho_n$ has to be
idempotent, and $\mu=\rho_n\otimes \delta_n$ is of the desired form.

Conversely, any probability of the form $\mu=\rho_n\otimes\delta_n$ with
$n\in\mathbb{Z}_+\cup\{\infty\}$ and $\rho_n$ idempotent is idempotent.
\end{proof}

For $k\in\mathbb{Z}$ and $n\in\mathbb{Z}_+$, define functions $\Theta^k_n:M\to
\mathbb{C}$ by
\[
\Theta^k_n(z,m) = \left\{\begin{array}{cll}
z^k & \mbox{ if }& m\ge n, \\
0  & \mbox{ if }& m < n.
\end{array}\right.
\]
The span of these functions is dense in $C(M)$, and they satisfy
\[
\Theta^k_n\Theta^\ell_m=\Theta^{k+\ell}_{\max(n,m)}, \quad (\Theta^k_n)^*=\Theta^{-k}_n, \quad\mbox{ and }\quad
\varepsilon(\Theta^k_n)=\Theta^k_n(e_M)=1.
\]
For their coproduct, we have
\begin{eqnarray*}
\Delta \Theta^k_n \big((z_1,m_1),(z_2,m_2)\big) &=&
\Theta_n^k\big(z_1z_2,\min(m_1,m_2)\big) \\
&=& \left\{\begin{array}{cll}
(z_1z_2)^k & \mbox{ if } & m_1,m_2\ge n, \\
0 & \mbox{ else},&
\end{array}\right.
\end{eqnarray*}
i.e.\ $\Delta \Theta^k_n=\Theta^k_n\otimes\Theta^k_n$.

\begin{proposition}
The semigroup $M$ is a quantum quotient semigroup of $U_0(2)$, in the sense
that there exists an injective $*$-algebra homomorphism $j$ from ${\rm Pol}(M)
:={\rm span}\{\Theta^k_n;n\in\mathbb{Z}_+,k\in\mathbb{Z}\}$ to ${\rm Pol}(U_0(2))$ such that
\[
\Delta_M\circ j = (j\otimes j)\circ \Delta.
\]
\end{proposition}
\begin{proof}
For $n\in\mathbb{Z}_+$ and $k\in\mathbb{Z}$, define $E^k_n=(\alpha^*)^n\alpha^nv^k\in {\rm Pol}(U_0(2))$. From 
the defining relations of $U_0(2)$, one can check that the $E^k_n$ satisfy the same $*$-algebraic and coalgebraic 
relations as the $\Theta_n^k$, i.e.\ $j(\Theta^k_n)=E_n^k$ defines a $*$-bialgebra homomorphism $j:{\rm 
Pol}(M)\to {\rm Pol}(U_0(2))$.

Let us show that $j$ is injective. Assume there exists a non-zero function $f=\sum_{k,n}
\lambda_{k,n}\Theta^k_n$ such that $j(f)=0$. Let $n_0$ be the smallest integer
  for which there exists a $k\in\mathbb{Z}$ such that
  $\lambda_{k,n}\not=0$. Take the representation $\pi=\pi_0\otimes{\rm id}_{L^2(\mathbb{T})}$ of
  $C(U_0(2))\cong C(SU_0(2))\otimes C(U(1))$ (recall that $\pi_0$ was defined in Section 2.3.1). Since $j(f)=0$, the operator
  $\chi=\pi(j(f))=\sum\lambda_{k,n}\pi_0\left((\alpha^*)^n\alpha^n\right)v^k$ has to
    vanish. Apply $\chi$ to $e_{n_0}\otimes 1$. Since
    $\pi_0\left((\alpha^*)^n\alpha^n\right)e_{n_0}=0$ for $n>n_0$, we get $\chi
      (e_{n_0}\otimes 1) = \sum_k \lambda_{n_0,k} e_0\otimes
      v^k\in\eufrak{h}\otimes L^2(\mathbb{T})$, which implies
      $\lambda_{n_0,k}=0$ for all $k$, in contradiction to the choice of $n_0$. 
\end{proof}

We can now give a description of all idempotent states on $U_0(2)$.
\begin{theorem}\label{theo-u02}
The following gives a complete list of the idempotent states on $U_0(2)$.
\begin{enumerate}
\item The idempotent states induced by the Haar measures on the two-dimensional torus $\mathbb{T}^2$ and its 
closed subgroups. If $\rho$ denotes the Haar measure of $\mathbb{T}^2$ or one of its closed subgroups, then the 
corresponding idempotent state $\phi_\rho$ is given by
\[
\phi_\rho\left((a^*)^rc^ka^sv^\ell\right) = \delta_{0k}\int_{\mathbb{T}^2} w_1^{s-r}
  w_2^\ell {\rm d}\rho(w_1,w_2)
\]
for $n,m\in\mathbb{Z}_+$, $k,\ell\in\mathbb{Z}$. This includes the counit of $U_0(2)$, for the trivial subgroup 
$\{1\}$ of $\mathbb{T}^2$. \item The family $\Psi_{n,m}=\psi_n\otimes \phi_m$, $n\in\mathbb{Z}_+$, 
$m\in\mathbb{N}\cup\{\infty\}$. Here $\phi_m$ is an idempotent state on $C(U(1))$, namely the Haar measure on 
$U(1)$ for $m=\infty$ and the idempotent state induced by the Haar measure of $\mathbb{Z}_m$ for 
$m\in\mathbb{N}$. And $\psi_n$, $n\in\mathbb{Z}_+$ is the idempotent state on $SU_0(2)$ defined by
\[
\psi_n\left((\alpha^*)^r\gamma^k\alpha^s\right) = \left\{\begin{array}{cll}
1 & \mbox{ if } & r=s\le n \mbox{ and }k=0 \\
0 & \mbox{ else}.&
\end{array}\right.
\]
\end{enumerate}
\end{theorem}

\begin{proof} 
Let $\phi:C(U_0(2))\to\mathbb{C}$ be an idempotent state on $U_0(2)$. Then
$\phi$ induces an idempotent state $\phi\circ j$ on $C(M)$. By Lemma
\ref{id-M}, $\phi\circ j$ is integration against a probability measure of the form $\rho\otimes \delta_{n}$ with
$n\in\mathbb{Z}_+\cup\{\infty\}$ and $\rho$ an idempotent measure on
$U(1)$. This determines $\phi$ on the subalgebra generated by $v$ and
$(a^*)^ra^r$, $r\in\mathbb{N}$: we have
\[
\phi\left((a^*)^ra^rv^k\right)=
\left\{
\begin{array}{lcl}
\rho(v^k) & \mbox{ if } & r\le n, \\
0 & \mbox{ else,}&
\end{array}
\right.
\]
for $k\in\mathbb{Z}$, $r\in\mathbb{Z}_+$.

\paragraph{Case (i): $n=\infty$}
For $k>0$ and any $r,s\in\mathbb{Z}_+$, $\ell\in\mathbb{Z}$, we have
\begin{eqnarray*}
\left|\phi\left((a^*)^r c^k a^s v^\ell\right)\right|^2 &\le&
\phi\left((a^*)^ra^r\right)\phi\left((a^*)^s (c^*)^kc^k a^s\right) \\
&=& \phi\left((a^*)^s a^s-(a^*)^{s+1} a^{s+1}\right) =0,
\end{eqnarray*}
and therefore $\phi$ vanishes on the $*$-ideal $\mathcal{I}_c$ generated by
$c$. As in the proof of Lemma \ref{lem-torus}, it follows that $\phi$ is
induced by an idempotent state on ${\rm Pol}(U_0(2))/\mathcal{I}_c\cong {\rm
  Pol}(\mathbb{T}^2)$, i.e.\ $\phi$ is of the form given in (1).
\paragraph{Case (ii): $n=0$}
Using again the Cauchy-Schwarz inequality, we get
\begin{eqnarray}
\left|\phi\left((a^*)^r c^k a^0 v^\ell\right)\right|^2 &\le&
\phi\left((a^*)^ra^r\right)\phi\left((a^*)^s a^s-(a^*)^{s+1}
  a^{s+1}\right),
\label{ineq-1} \\
\left|\phi\left((a^*)^r c^k a^s v^\ell\right)\right|^2 &\le& \phi\left((a^*)^r a^r-(a^*)^{r+1}a^{r+1}\right)\phi\left((a^*)^s a^s\right),
\label{ineq-2} \\
\left|\phi\left((a^*)^r a^s v^\ell\right)\right|^2 &\le&
\phi\left((a^*)^ra^r\right)\phi\left((a^*)^s a^s\right),
\label{ineq-3}
\end{eqnarray}
for $k,r,s\in\mathbb{Z}_+$, $\ell\in\mathbb{Z}$. This shows that
$\phi\left((a^*)^r c^k a^r v^\ell\right)$ vanishes, unless $r=s=0$.
different from $0$. But then $\phi=\phi\star\phi$ implies
\[
\phi(c^kv^\ell) = (\phi\otimes\phi)\big(\Delta(c^k v^\ell)\big) = \sum_{\kappa=0}^k \phi\left((a^*)^{k-\kappa}c^\kappa
  v^\ell\right)\phi\left(c^{k-\kappa}a^\kappa v^\ell\right) = 0
\]
for $k>0$, $s\in\mathbb{Z}$. By hermitianity
$\phi\left((c^*)^kv^\ell\right)=0$ for $k>0$, and $\phi$ has the form given in
(2) with $n=0$.
\paragraph{Case (iii): $n\in\mathbb{N}$}
We use once more the Cauchy-Schwarz inequality. For $k\not=0$, \eqref{ineq-1} and \eqref{ineq-2}
imply that $\phi\left((a^*)^r c^k a^s v^\ell\right)$ vanishes unless $r=s=n$. But then we can show
that $\phi=\phi\star\phi$ implies $\phi\left((a^*)^n c^k a^n v^\ell\right)=0$ in the same way as in
the previous case.

For $k=0$, we see from \eqref{ineq-3} that $\phi\left((a^*)^r a^s
  v^\ell\right)$ vanishes unless $r,s\le n$. The elements $(a^*)^r a^s
  v^\ell$ are group-like, therefore $\phi\left((a^*)^r a^s
  v^\ell\right)\in\{0,1\}$. If we can show $\phi\left((a^*)^r
  a^s v^\ell\right)\not=1$ for $r\not=s$, we are done, since then
$\phi\left((a^*)^r c^k a^s v^\ell\right)$ is non-zero only if $r=s\le n$.
We get $\phi\left((a^*)^r a^s v^\ell\right)= \delta_{rs}\rho(v^\ell)$ for
$r,s\le n$, i.e.\
$\phi$ has the form given in (2).

We show $\phi\left((a^*)^r
  a^s v^\ell\right)\not=1$ for $r\not=s$ by contradiction. Assume
there exists a triple $(r_0,s_0,\ell_0)$ such that $\phi\left((a^*)^r
  a^s v^\ell\right)=1$ and choose such a triple with maximal $r_0$. Set
\[
b=(a^*)^{r_0}a^{s_0}v^{\ell_0} + (a^*)^{s_0}a^{r_0}v^{-\ell_0} -1.
\]
Maximality of $r_0$ implies $r_0>s_0$ and
\[
\phi\left((a^*)^{2r_0-s_0}a^{s_0}v^{2\ell_0}\right)=\phi\left((a^*)^{s_0}a^{2r_0-s_0}v^{-2\ell_0}\right)=0,
\]
therefore we get
\begin{eqnarray*}
\phi(b^*b) &=&
\phi\left((a^*)^{r_0}a^{r_0}\right)
+\phi\left((a^*)^{s_0}a^{s_0}\right)
+1 \\
&& +\phi\left((a^*)^{2r_0-s_0}a^{s_0}v^{2\ell}\right)
+\phi\left((a^*)^{s_0}a^{2r_0-s_0}v^{-2\ell}\right) \\
&& -2\phi\left((a^*)^{r_0}a^{s_0}v^{\ell}\right)
-2\phi\left((a^*)^{s_0}a^{r_0}v^{-\ell}\right) \\
&=& -1,
\end{eqnarray*}
which is clearly a contradiction to the positivity of $\phi$.

Conversely, using the formulas
\begin{eqnarray*}\Delta\big((a^*)^r c^k a^s v^\ell\big) &=& \sum_{\kappa=0}^k
(a^*)^{r+k-\kappa}c^\kappa a^s v^\ell \otimes (a^*)^r
c^{k-\kappa}a^{s+\kappa}v^\ell, \\
\Delta\big((a^*)^r (c^*)^k a^s v^\ell\big) &=& \sum_{\kappa=0}^k
(a^*)^r (c^*)^\kappa a^{s+k-\kappa} v^\ell \otimes (a^*)^{r+\kappa}
(c^*)^{k-\kappa}a^s v^\ell,
\end{eqnarray*}
for $r,s,k\in\mathbb{Z}_+$, $\ell\in\mathbb{Z}$,
for the coproduct in ${\rm Pol}(U_0(2))$, one can check that all states given
in the theorem are indeed idempotent.
\end{proof}

\begin{remark}
The state $\Psi_{0,\infty}$ introduced in the theorem above can be considered as the Haar state on $U_0(2)$
since it is invariant, i.e.\
\[
\Psi_{0,\infty}\star f = f\star \Psi_{0,\infty} = f(1) \Psi_{0,\infty}
\]
for any $f\in C(U_0(2))^*$. But $\Psi_{0,\infty}$ is not faithful. Its algebraic null space
\begin{eqnarray*}
N_{\Psi_{0,\infty}} &=& \{u\in {\rm
  Pol}(U_0(2));\Psi_{0,\infty}(u^*u)=0\} \\
&=& {\rm span}\,\{(a^*)^kc^m a^\ell v^n;
k\in\mathbb{Z}_+,m,n\in\mathbb{Z},\ell\ge 1\}
\end{eqnarray*}
is a left ideal, but not self-adjoint or two-sided. It is a subcoalgebra,
i.e.\ we have
\[
\Delta N_{\Psi_{0,\infty}} \subseteq
N_{\Psi_{0,\infty}}\odot N_{\Psi_{0,\infty}},
\]
but it is not a coideal, since the counit does not vanish on $\mathcal{N}_{\Psi_{0,\infty}}$.
\end{remark}

The complete description of the
idempotent states on $SU_0(2)$ follows now directly from Theorem \ref{theo-u02} and the comments before Theorem \ref{theo-uq2}.

\begin{theorem}\label{theo-su02}
The following gives a complete list of the idempotent states on $SU_0(2)$.
\begin{enumerate}
\item
The family $\phi_n$, $n\in\mathbb{N}\cup\{\infty\}$ where $\phi_1$ is the
counit, $\phi_\infty$ the idempotent induced by the Haar state
of quantum subgroup $U(1)$, and $\phi_n$, $2\le n < \infty$, denotes the
idempotent state induced by the Haar state on the quantum subgroup $\mathbb{Z}_n$.
\item
The family $\psi_n$, $n\in\mathbb{Z}_+$ defined by
\[
\psi_n\big((\alpha^*)^r \gamma^k \alpha^s\big) =
\left\{\begin{array}{lll}
1 & \mbox{ if } & k=0, r=s\le n, \\
0 & \mbox{ else.}&
\end{array}\right.
\]
for $r,s\in\mathbb{Z}_+$, $k\in\mathbb{Z}$, with the convention $\gamma^{-k}=(\gamma^*)^k$.
\end{enumerate}
\end{theorem}

Similarly using the conditional expectation introduced in Paragraph
\ref{subsub-cond-exp}, we can derive a complete classification of the
idempotent states on $SO_0(3)$. The proof is identical to the proof of Theorem
\ref{theo-soq3} and therefore omitted.

\begin{theorem}\label{theo-so03}
All idempotent states on $SO_0(3)$ arise as restrictions of idempotent states
on $SU_0(2)$. Moreover $\phi_n|_{C(SO_0(3))} =
\phi_{2n}|_{C(SO_0(3))}$ for $n$ an odd integer.
\end{theorem}

\section{The idempotent states on $SU_q(2)$ as elements of the dual and associated quantum hypergroups}
\label{hypergr} Let $\QG$ be a compact quantum group and $\mathcal{A}=\Pol$ the corresponding Hopf $*$-algebra 
dense in $C(\QG)$. Then $\mathcal{A}$ is (the algebra of functions on) an \emph{algebraic quantum group} in the 
sense of Van Daele (\cite{Van Daele 98}), and so is its dual $\hat{\mathcal{A}}$, given by the functionals of the 
form $h(\cdot\,a)$ with $a\in\mathcal{A}$. By \cite[Lemma 3.1]{franz+skalski08a}, an idempotent state $\phi$ on 
$C(\QG)$ defines a group-like projection $p_\phi$ in the multiplier algebra $M(\hat{\mathcal{A}})$ of the dual, 
and therefore, by \cite[Theorem 2.7]{landstad+vandaele07} and \cite[Theorem 2.4]{franz+skalski08a} an 
\emph{algebraic quantum hypergroup} $\hat{\mathcal{A}}_{p_\phi}$. As an algebra, $\hat{\mathcal{A}}_{p_\phi}= 
p_\phi\hat{\mathcal{A}}p_\phi$, and the coproduct of $\hat{\mathcal{A}}_{p_\phi}$ is given by
\[
\hat{\Delta}_{p_\phi} = (p_\phi\otimes p_\phi)\hat{\Delta}(a)(p_\phi\otimes p_\phi)
\]
for $a\in \hat{\mathcal{A}}_{p_\phi}$, where $\hat{\Delta}$ denotes the coproduct of
$\hat{\mathcal{A}}$.

Let $q\in (-1,0)\cup(0,1)$. In
this section we will consider the case of the compact quantum group $SU_q(2)$
and describe the algebraic quantum hypergroups associated to its idempotent
states. Note that in this case the dense Hopf $*$-algebra is
$\mathcal{A}={\rm Pol}(SU_q(2))={\rm span}\,\{u^{(s)}_{k\ell}: s\in\frac{1}{2}\mathbb{Z}_+,-s\le
    k,\ell\le s\}$. We will use the basis
\[
e^{(s)}_{kl} = \frac{1-q^{2(2s+1)}}{q^{2(s-k)}(1-q^2)}h\big((u^{(s)}_{k\ell})^*\,\cdot\,\big)
\]
for $\hat{\mathcal{A}}$, which can be thought of as the algebra of trigonometric polynomials on the quantum group dual $\widehat{SU_q(2)}$. Using the orthogonality
relation
\[
h\left((u^{(s)}_{k\ell})^* u^{(s')}_{k'\ell'}\right) = \delta_{ss'} \delta_{kk'}\delta_{\ell\ell'}\frac{q^{2(s-k)}(1-q^2)}{1-q^{2(2s+1)}},
\]
for $s,s'\in\frac{1}{2}\mathbb{Z}_+$, $-s\le k,\ell\le s$, $-s'\le k',\ell'\le
s'$, cf.\ \cite[Eq.\ (5.12)]{koornwinder89},
we can check that this basis is dual to the basis $\{u^{(s)}_{k\ell}:
s\in\frac{1}{2}\mathbb{Z}_+,-s\le k,\ell\le s\}$ of ${\rm Pol}(SU_q(2))$. The
algebraic quantum group
$\hat{\mathcal{A}}={\rm Pol}(\widehat{SU_q(2)})$ is of discrete type and equal to the algebraic direct sum
\[
\hat{\mathcal{A}} =\bigoplus_{s\in \frac{1}{2}\mathbb{Z}_+} M_{2s+1}.
\]
The $e^{(s)}_{k\ell}$ form a basis of matrix units for $M_{2s+1}={\rm
  span}\,\{e^{(s)}_{kl}: -s\le k,\ell\le s\}$.

The Haar state $h$ and the counit $\varepsilon$ give the elements $p_h=1$ and
$p_\varepsilon=e^{(0)}_{00}$ in $M(\hat{\mathcal{A}})$, and the associated
algebraic quantum hypergroups are $\hat{\mathcal{A}}_{p_h}=\hat{\mathcal{A}}$ and $\hat{\mathcal{A}}_{p_\varepsilon}=\mathbb{C}$.

The remaining cases are more interesting.

\subsection{The idempotent state $\phi_2 = h_{\mathbb{Z}_2}\circ j$ induced
  by the quantum  subgroup $\mathbb{Z}_2$}

We have
\[
p_{\phi_2} = \sum_{s=0}^\infty \sum_{k=-s}^s e^{(s)}_{kk} = \sum_{s=0}^\infty
1_{2s+1},
\]
i.e.\ $p_{\phi_2}$ is the sum of the identity matrices from the
odd-dimensional matrix algebras $M_{2s+1}$, $s\in\mathbb{Z}_+$. This
projection is in the center of $M(\hat{\mathcal{A}})$, therefore $\hat{\mathcal{A}}_2=p_{\phi_2}\hat{\mathcal{A}}p_{\phi_2}$
will be an algebraic quantum group. We get
\[
\hat{\mathcal{A}}_2= \bigoplus_{s\in \mathbb{Z}_+} M_{2s+1} = {\rm Pol}\,
\widehat{SO_q(3)},
\]
i.e.\ $\hat{\mathcal{A}}_{p_{\phi_2}}$ is the algebra of trigonometric polynomials on the algebraic quantum group dual to $SO_q(3)$. This is to be expected as
$\mathbb{Z}_2$ is the only nontrivial normal quantum subgroup of $SU_q(2)$ and $SO_q(3)$ is the corresponding quotient quantum group.

\subsection{The idempotent state $\phi_\infty = h_{U(1)}\circ j$ induced by the quantum subgroup  $U(1)$}

Here
\[
p_{\phi_{\infty}}= \sum_{s=0}^\infty e^{(s)}_{00}
\]
and this projection is not central. We get
\[
\hat{\mathcal{A}}_{p_\infty} = \bigoplus_{s=0}^\infty \mathbb{C}
\]
which is a commutative algebraic quantum hypergroup of discrete type. This is the dual of the hypergroup introduced in \cite[Section 7]{koornwinder91}.

\subsection{The idempotent states $\phi_n = h_{\mathbb{Z}_{n}}\circ j$, $3\le n<\infty$}

The remaining cases also give non-central projections,
\begin{eqnarray*}
p_{\phi_{2n}} &=& \sum_{s=0}^\infty
\,\sum_{k=-\lceil\frac{s}{n}\rceil}^{\lceil\frac{s}{n}\rceil} e^{(s)}_{nk,nk},\\
p_{\phi_{2n+1}} &=&
\sum_{s\in\frac{1}{2}\mathbb{Z}_+}\,\sum_{\genfrac{.}{.}{0pt}{}{-s\le k \le
    s}{2k\equiv 0\!\!\!\!\mod 2n+1}} e^{(s)}_{kk},
\end{eqnarray*}
for $1\le n<\infty$, cf.\ Equations \eqref{z2n-idemp} and \eqref{z2n+1-idemp}. They lead to non-commutative algebraic
quantum hypergroups $\hat{\mathcal{A}}_{n}=p_{\phi_n}\hat{\mathcal{A}}p_{\phi_n}$,
\begin{eqnarray*}
\hat{\mathcal{A}}_{2n} &=& \bigoplus_{k=0}^\infty  n M_{2k+1}, \\
\hat{\mathcal{A}}_{2n+1} &=& \bigoplus_{k=0}^\infty (2n+1) M_k, \\
\end{eqnarray*}
of discrete type (in the formulas above $nM_{2k+1}$ denotes $n$ direct copies of the matrix algebra $M_{2k+1}$ and similarly
$(2n+1)M_k$ denotes $2n+1$ direct copies of the matrix algebra $M_{k}$).

Note that as the quantum subgroups consider in the last two paragraphs are not normal, the objects we obtain have only the
quantum hypergroup structure (and can be informally thought of as duals of quantum hypergroups obtained via the double coset
construction, \cite{chapovsky+vainerman99}).

\renewcommand{\theequation}{\Alph{section}.\arabic{equation}}

\renewcommand{\thesection}{\Alph{section}}

\setcounter{section}{0}

\renewcommand{\thetheorem}{\Alph{section}.\arabic{theorem}}

\section*{Appendix}

\setcounter{section}{1} \setcounter{theorem}{0}

The goal of the appendix is to provide an alternative proof of coamenability of the deformations of classical 
compact Lie groups. To facilitate the discussion, for a compact quantum group $QG$ we use $\widehat{\mathbb G}$ 
to denote the dual of $\mathbb{G}$ (which is a discrete quantum group), $C(\mathbb{G})_{\rm red}$ and 
$C(\mathbb{G})$ to denote respectively the reduced and universal $C^*$-algebras associated with $\mathbb{G}$ and 
$L^{\infty}(\mathbb{G})$ to denote the corresponding von Neumann algebra (we refer for example to 
\cite{tomatsu07} for precise definitions).  Note that contrary to the main body of the paper we do not assume 
that the Haar state on $\mathbb{G}$ is faithful, so that $\mathbb{G}$ need not be in the reduced form. We adopt 
the following definition (\cite{bedos+murphy+tuset02}, \cite{tomatsu06}).

\begin{definition} \label{defamen}
A compact quantum group $\mathbb{G}$ is said to be \emph{coamenable} if the dual quantum group $\widehat{\mathbb G}$ is amenable,
that is, $L^\infty(\widehat{\mathbb G})$ has an invariant mean.
\end{definition}

The following result gives a useful criterion to check coamenability:

\begin{theorem}\cite[Theorem 4.7]{bedos+murphy+tuset02}, \cite[Corollary 3.7, Theorem 3.8]{tomatsu06}
A compact quantum group $\mathbb{G}$ is coamenable if and only if there exists a counit on
$C(\mathbb{G})_{\rm red}$ if and only if there exists a $*$-homorphism from $C(\mathbb{G})_{\rm
red}$ onto $\mathbb{C}$. \label{coamequiv}\end{theorem}

The second equivalence is fairly easy to show, in the first the forward implication was established
in \cite{bedos+murphy+tuset02} and the backward implication in \cite{tomatsu06}.

Let $\mathbb{G}$ be a classical compact Lie group and $\mathbb{G}_q$ the $q$-deformation with the
parameter $-1<q<1$, $q\neq0$ (see \cite{korogodski+soibelman98}). The function algebra
$C(\mathbb{G}_q)$ is the universal $C^*$-algebra generated by certain polynomial elements. The Haar
state is denoted by $h$.

The following theorem was proved by T. Banica (\cite[Corollary 6.2]{Banica99}). We present another proof using L.I.\
Korogodski-Y.S.\ Soibelman's results on the representation theory of $C(\mathbb{G}_q)$.

\begin{theorem}
The quantum group $\mathbb{G}_q$ is coamenable.
\end{theorem}
\begin{proof}
Let us introduce the left ideal $N_h:=\{a\in C(\mathbb{G}_q)\mid h(a^*a)=0\}$, which is in fact an ideal of 
$C(\mathbb{G}_q)$. The reduced version  $C(\mathbb{G}_q)_{\rm red}$ is defined as the quotient 
$C(\mathbb{G}_q)/N_h$. By Theorem \ref{coamequiv} to show that $\mathbb{G}_q$ is coamenable it suffices to show 
that the $C^*$-algebra $C(\mathbb{G}_q)_{\rm red}$ has a character.

Consider an irreducible representation $\pi\colon C(\mathbb{G}_q)_{\rm red}\to B(H_\pi)$. Composing
this map with the canonical surjection $\rho\colon C(\mathbb{G}_q)\to C(\mathbb{G}_q)_{\rm red}$,
we get an irreducible representation $\pi\circ\rho$ of $C(\mathbb{G}_q)$. Thanks to \cite[Theorem
6.2.7 (3), \S3]{korogodski+soibelman98}, we may assume that $\pi\circ\rho$ is of the following
form:
\[
\pi\circ\rho = (\pi_{s_{i_1}}\otimes\cdots\otimes\pi_{s_{i_k}}\otimes\pi_t)\circ\Delta^{(k)}
\]
or $\pi\circ\rho=\pi_t$, where $s_{i_1}\cdots s_{i_k}$ is the reduced decomposition in the Weyl group of $\mathbb{G}$, and $t\in
T$, the maximal torus of $\mathbb{G}$. In the latter case $\pi$ is a one-dimensional representation. In the former case, we
remark that the counit of $C(\mathbb{G}_q)$ factors through $\textup{Im}\, \pi_{s_i}$ for every $i$, that is, there exists
$\eta_i\colon \textup{Im}\, \pi_{s_i}\rightarrow \mathbb{C}$ such that $\eta_i\circ\pi_{s_i}=\varepsilon$ (See the argument in
\cite[p. 294]{tomatsu07}).

Then we introduce a representation
$\tilde{\pi}:=(\eta_{i_1}\otimes\cdots\otimes\eta_{i_k}\otimes\textup{id})\circ\pi$ of
$C(\mathbb{G}_q)_{\rm red}$, which is well-defined and one-dimensional. Indeed,
\[
\tilde{\pi}\circ \rho =(\varepsilon\otimes\cdots\otimes\varepsilon\otimes\pi_t)\circ\Delta^{(k)} =\pi_t.
\]
Thus we have proved in each case the existence of a one-dimensional representation of the $C^*$-algebra $C(\mathbb{G}_q)_{\rm
red}$, and $\mathbb{G}_q$ is coamenable.
\end{proof}

\section*{Acknowledgements}

Big part of this research was done at the Mathematisches Forschungsinstitut Oberwolfach
during a  stay within the Research in Pairs Programme from September 7 to September 20, 2008. We thank the MFO for giving us this opportunity.

This work was started while U.F.\ was visiting the Graduate School of
Information Sciences of Tohoku University as Marie-Curie fellow. He would like
to express his gratitude to Professors Nobuaki Obata, Fumio Hiai, and the
other members of the GSIS for their hospitality.


\newcommand{\etalchar}[1]{$^{#1}$}

\end{document}